\documentclass[10pt]{amsart}

\usepackage{spectralsequences}

\usepackage{mathtools}
\usepackage{amsthm}
\usepackage{mathrsfs} 
\usepackage{amssymb}
\usepackage[all]{xy}
\usepackage[nobottomtitles*]{titlesec}
\usepackage{titletoc}
\usepackage{bbm}
\usepackage{upgreek}
\usepackage{ marvosym}
\usepackage{multicol}
\usepackage{rotating}
\usepackage{pdflscape}
\usepackage{ifluatex}
\usepackage[pdfencoding=auto, psdextra]{hyperref}
\usepackage[noabbrev,nameinlink]{cleveref}
\usepackage{enumerate}
\usepackage{soul}
\usepackage{graphicx}
\usepackage{stmaryrd}
\usepackage{todonotes}

\usepackage{verbatim}
\usepackage{tikz-cd}
\makeatletter
\def\@makefntext{\noindent\@makefnmark}
\tikzcdset{
    cong/.style={"\cong" {sloped, description, yshift=0pt,#1}, phantom},
    simeq/.style={"\simeq" {sloped, description, yshift=0pt,#1}, phantom},
    snake/.style={
        out= east, in=west,
        to path={
            \pgfextra{
                \pgfextractx{\pgf@xa}{\pgfpointanchor{\tikztostart}{east}}
                \pgfextractx{\pgf@xb}{\pgfpointanchor{\tikztotarget}{west}}
                \pgfextracty{\pgf@ya}{\pgfpointanchor{\tikztostart}{center}}
                \pgfextracty{\pgf@yb}{\pgfpointanchor{\tikztotarget}{center}}
                \edef\tikzstartx{\the\pgf@xa}
                \edef\tikzendx{\the\pgf@xb}
                \edef\midy{\the\dimexpr0.5\dimexpr\pgf@ya\relax +0.5\dimexpr\pgf@yb\relax}
            }
            to[in=0,out=180,looseness=0.5] (\tikzstartx,\midy)
            -| ([xshift=-2ex]\tikztotarget.west)
            -- (\tikztotarget)}
    }
}
\makeatother
\usepackage{sseq}
\usepackage{tikzcdintertext}
\usepackage{etoolbox}
\usepackage{xpatch}

\crefname{diagram}{diagram}{diagrams}
\creflabelformat{diagram}{#2(#1)#3}
\crefname{sseq}{}{}
\crefname{equation}{}{}
\Crefname{equation}{Equation}{Equations}
\crefname{equation-b}{Equation}{Equations}
\creflabelformat{sseq}{#2(#1)#3}
\usepackage[normalem]{ulem}
\newcommand{\stkout}[1]{\ifmmode\text{\sout{\ensuremath{#1}}}\else\sout{#1}\fi}

\makeatletter

\pretocmd{\maketitle}{%
    \edef\@title{\unexpanded{\protect\Large}\unexpanded\expandafter{\@title}}%
    \edef\authors{\unexpanded{\protect\normalsize}\unexpanded\expandafter{\authors}}%
}{}{\error}

\def\addlabeltolink#1{\addlabeltolink@#1}
\def\addlabeltolink@#1#2#3#4{#1{#2}{#3}{\thecontentslabel. #4}}

\dottedcontents{section}[0em]{}{2.3em}{1pc}
\dottedcontents{section}[3.8em]{}{2.3em}{1pc}

\titleformat{\section}[block]{\centering\bfseries\Large}{\thetitle. }{0pt}{}

\def\@secnumpunct{. }
\xpatchcmd{\proof}{\topsep6\p@\@plus6\p@\relax}{\topsep0pt\relax}{}{\error}
\makeatother
\usepackage[parfill]{parskip}

\allowdisplaybreaks[1]

\newcommand{\ull}[1]{\underline{#1}}

\newcommand{\mr}[1]{\mathrm{#1}}
\usepackage[nobottomtitles*]{titlesec}
\usepackage{titletoc}

\renewcommand{\psi}{\uppsi}

\newcommand{\PP}{\mathbb{P}}

\newcommand{\tmf}{\mathsf{tmf}}

\newcommand{\KO}{\mr{KO}}

\newcommand{\RO}{\mr{RO}}

\newcommand*{\rom}[1]{\expandafter\@slowromancap\romannumeral #1@} \makeatother

\hypersetup{%
  bookmarksnumbered=true,%
  bookmarks=true,%
  colorlinks=true,%
  pdfnewwindow=true,%
  pdfstartview=FitBH%
}

\hypersetup{
    colorlinks,
    linkcolor={DarkBrown},
    citecolor={DarkBlue},
    urlcolor={blue!80!black}
}

\def\@url#1{{\tt\def~{\lower3.5pt\hbox{\char'176}}\def\_{\char'137}#1}}

\newtheorem{main}{Main Theorem}
\crefname{main}{main theorem}{main theorems}
\Crefname{main}{Main Theorem}{Main Theorems}

\newtheorem{thm}[equation]{Theorem}
\newtheorem{cor}[equation]{Corollary}
\newtheorem{lem}[equation]{Lemma}

\newtheorem{prop}[equation]{Proposition}

\theoremstyle{definition}

\newtheorem{rmk}[equation]{Remark}

\newtheorem{notn}[equation]{Notation}

\crefname{thm}{theorem}{theorems}
\Crefname{thm}{Theorem}{Theorems}
\crefname{cor}{corollary}{corollaries}
\Crefname{cor}{Corollary}{Corollaries}
\crefname{lem}{lemma}{lemmas}
\Crefname{lem}{Lemma}{Lemmas}
\crefname{prop}{proposition}{propositions}
\Crefname{prop}{Proposition}{Propositions}
\crefname{defn}{definition}{definitions}
\Crefname{defn}{Definition}{Definitions}
\crefname{ex}{example}{examples}
\Crefname{ex}{Example}{Examples}
\crefname{exs}{example}{examples}
\Crefname{exs}{Example}{Examples}
\crefname{rmk}{remark}{remarks}
\Crefname{rmk}{Remark}{Remarks}
\crefname{conj}{conjecture}{conjectures}
\Crefname{conj}{Conjecture}{Conjectures}
\crefname{obs}{observation}{observations}
\Crefname{obs}{Observation}{Observations}
\crefname{claim}{claim}{claims}
\Crefname{claim}{Claim}{Claims}
\crefname{prob}{problem}{problems}
\Crefname{prob}{Problem}{Problems}
\crefname{warn}{warning}{warnings}
\Crefname{warn}{Warning}{Warnings}
\crefname{notn}{notation}{notations}
\Crefname{notn}{Notation}{Notations}
\crefname{formula}{formula}{formulas}
\Crefname{formula}{Formula}{Formulas}
\crefname{construction}{construction}{constructions}
\Crefname{construction}{Construction}{Constructions}
\crefname{question}{question}{questions}
\Crefname{question}{Question}{Questions}
\crefname{assump}{assumption}{assumptions}
\Crefname{assump}{Assumption}{Assumptions}

\AddToHook{env/thm/begin}{\crefalias{equation}{thm}}
\AddToHook{env/cor/begin}{\crefalias{equation}{cor}}
\AddToHook{env/lem/begin}{\crefalias{equation}{lem}}
\AddToHook{env/prob/begin}{\crefalias{equation}{prob}}
\AddToHook{env/prop/begin}{\crefalias{equation}{prop}}
\AddToHook{env/obs/begin}{\crefalias{equation}{obs}}
\AddToHook{env/conj/begin}{\crefalias{equation}{conj}}
\AddToHook{env/claim/begin}{\crefalias{equation}{claim}}
\AddToHook{env/defn/begin}{\crefalias{equation}{defn}}
\AddToHook{env/formula/begin}{\crefalias{equation}{formula}}
\AddToHook{env/ex/begin}{\crefalias{equation}{ex}}
\AddToHook{env/exs/begin}{\crefalias{equation}{exs}}
\AddToHook{env/rmk/begin}{\crefalias{equation}{rmk}}
\AddToHook{env/warn/begin}{\crefalias{equation}{warn}}
\AddToHook{env/notn/begin}{\crefalias{equation}{notn}}
\AddToHook{env/construction/begin}{\crefalias{equation}{construction}}
\AddToHook{env/question/begin}{\crefalias{equation}{question}}

\newtheorem*{thm*}{Theorem}
\newtheorem*{cor*}{Corollary}
\newtheorem*{lem*}{Lemma}
\newtheorem*{prop*}{Proposition}
\newtheorem*{not*}{Notation}
\newtheorem*{guess*}{Guess}

\newtheorem*{defn*}{Definition}
\newtheorem*{ex*}{Example}
\newtheorem*{exs*}{Examples}
\newtheorem*{rmk*}{Remark}
\newtheorem*{claim*}{Claim}
\newtheorem*{exer*}{Exercise}

\numberwithin{equation}{section}
\numberwithin{figure}{section}

\usepackage[T1]{fontenc}
\makeatletter
\newif\ifkp@upRm
\DeclareSymbolFont{Letters}{OML}{jkp}{m}{n}
\DeclareMathSymbol{\partialup}{\mathord}{Letters}{128}

\DeclareMathOperator{\Sp}{\mathcal{S}p}

\DeclareMathOperator{\Top}{\mathcal{T}op}

\DeclareMathOperator{\Th}{{\sf Th}} 

\DeclareMathOperator{\res}{{\sf res}}

\def\makecommands#1#2#3{
    \bgroup
    \def\tempcmdname##1{#1}
    \def\tempcmdbody##1{#2}
    \def\\##1{\expandafter\xdef\csname\tempcmdname{##1}\endcsname{\unexpanded\expandafter{\tempcmdbody{##1}}}}
    #3
    \egroup
}

\def\upperalphabet{\\A\\B\\C\\D\\E\\F\\G\\H\\I\\J\\K\\L\\M\\N\\O\\P\\Q\\R\\S\\T\\U\\V\\W\\X\\Y\\Z}
\def\loweralphabet{\\a\\b\\c\\d\\e\\f\\g\\h\\i\\j\\k\\l\\m\\n\\o\\p\\q\\r\\s\\t\\u\\v\\w\\x\\y\\z}
\def\lowergreekalphabet{\\\alpha\\\beta\\\gamma\\\delta\\\epsilon\\\zeta\\\eta\\\theta\\\kappa\\\lambda\\\mu\\\nu
    \\\xi\\\pi\\\rho\\\sigma\\\tau\\\upsilon\\\psi\\\chi\\\phi\\\omega}

\makecommands{#1mr}{\mathrm{#1}}{\upperalphabet}  
\makecommands{#1cal}{\mathcal{#1}}{\upperalphabet}  
\makecommands{#1#1}{\mathbb{#1}}{\upperalphabet} 
\makecommands{#1bar}{\overline{#1}}{\upperalphabet}
\makecommands{#1scr}{\mathscr{#1}}{\upperalphabet}
\makecommands{#1twee}{\widetilde{#1}}{\upperalphabet\loweralphabet} 
\makecommands{sf#1}{{\sf #1}}{\upperalphabet\loweralphabet}
\makecommands{#1frak}{\mathfrak{#1}}{\upperalphabet}  
\makeatletter
\makecommands{\expandafter\@gobble\string#1bar}{\overline{#1}}{\lowergreekalphabet}
\makecommands{\expandafter\@gobble\string#1twee}{\widetilde{#1}}{\lowergreekalphabet}
\makeatother

\newcommand{\sma}{\wedge}

\definecolor{limegreen}{rgb}{0.2, 0.8, 0.2}
\definecolor{darkmagenta}{rgb}{0.55, 0.0, 0.55}
\definecolor{lavenderrose}{rgb}{0.91, 0.33, 0.5}
\definecolor{goldenpoppy}{rgb}{0.99, 0.76, 0.0}
\definecolor{seagreen}{rgb}{0.1, 0.4, 0.1}
\definecolor{maroon}{RGB}{128,0,0}
\definecolor{darkviolet}{RGB}{148,0,211}
\definecolor{darkbrown}{rgb}{.55,.3,.2}
\definecolor{snowblue}{rgb}{0.3, 0.4, 0.8}
\definecolor{DarkBlue}{rgb}{.1, 0.35, 0.6} 
\definecolor{DarkBrown}{rgb}{.5, 0.2, 0.2}

\makeatletter

\def\tikzcdequalsignoffset{0.1em}

\def\findedgesourcetarget#1#2{
    \let\sourcecoordinate\pgfutil@empty
    \ifx\tikzcd@startanchor\pgfutil@empty 
        \def\tempa{\pgfpointanchor{#1}{center}}
    \else
        \edef\tempa{\noexpand\pgfpointanchor{#1}{\expandafter\@gobble\tikzcd@startanchor}} 
        \let\sourcecoordinate\tempa
    \fi
    \ifx\tikzcd@endanchor\pgfutil@empty 
        \def\tempb{\pgfpointshapeborder{#2}{\tempa}}
    \else
        \edef\tempb{\noexpand\pgfpointanchor{#2}{\expandafter\@gobble\tikzcd@endanchor}}
    \fi
    \let\targetcoordinate\tempb
    \ifx\sourcecoordinate\pgfutil@empty%
        \def\sourcecoordinate{\pgfpointshapeborder{#1}{\tempb}}%
    \fi
}

\tikzset{/tikz/commutative diagrams/equal/.style=equals,
    /tikz/commutative diagrams/equals/.style = {
    -,
    to path={\pgfextra{
        \findedgesourcetarget{\tikzcd@ar@start}{\tikzcd@ar@target} 
        \ifx\tikzcd@startanchor\pgfutil@empty
            \def\tikzcd@startanchor{.center}
        \fi
        \ifx\tikzcd@endanchor\pgfutil@empty
            \def\tikzcd@endanchor{.center}
        \fi
        \pgfmathanglebetweenpoints{\pgfpointanchor{\tikzcd@ar@start}{\expandafter\@gobble\tikzcd@startanchor}}{\pgfpointanchor{\tikzcd@ar@target}{\expandafter\@gobble\tikzcd@endanchor}}
        \pgftransformrotate{\pgfmathresult}
        \pgfpathmoveto{\pgfpointadd{\sourcecoordinate}{\pgfpoint{0}{\tikzcdequalsignoffset}}}
        \pgfpathlineto{\pgfpointadd{\targetcoordinate}{\pgfpoint{0}{\tikzcdequalsignoffset}}}
        \pgfpathmoveto{\pgfpointadd{\sourcecoordinate}{\pgfpoint{0}{-\tikzcdequalsignoffset}}}
        \pgfpathlineto{\pgfpointadd{\targetcoordinate}{\pgfpoint{0}{-\tikzcdequalsignoffset}}}
        \pgfusepath{draw}
}}}}

\makeatother

\title{New $\mathrm{C}_2$-equivariant infinite families and the $\mathrm{H}\underline{\mathbb{F}}_2$-Hurewicz image}
\author{M. Adams}\address{New Mexico State University}\email{mtadams0@nmsu.edu}
\author{P. Bhattacharya}\address{New Mexico State University}\email{prasit@nmsu.edu}
\author{A. Bhattacherjee}\address{New Mexico State University}\email{arannya5@nmsu.edu}
\author{A. Waugh}\address{ New Mexico State University}\email{awaugh@nmsu.edu}

\thanks{This work is supported by the National Science Foundation under grants DMS-2305016 and DMS-2603448.}

\begin{document}

\begin{abstract} 
Let $a$ denote the stabilization of the Euler class of the sign representation of $\Cmr_2$. In this paper, we show that the Behrens-Shah periodic $u$-self-equivalence of $C(a^i)$ -- a direct consequence of J.F. Adams' calculation of James periodicity --  determines the  $\mathrm{H}\underline{\mathbb{F}}_2$-Hurewicz image of the $\Cmr_2$-equivariant sphere spectrum, recovering a recent theorem of Guo--Li--Lu--Wu--Xu--Yang--Zhang via a different method. Furthermore, we use the $\Cmr_2$-equivariant periodic $u$-self-equivalences on the cofiber spectra $C(a^i)$ to construct new infinite families in the $\Cmr_2$-equivariant stable stems whose $\mathrm{H}\underline{\mathbb{F}}_2$-Hurewicz image is trivial. \end{abstract}

\maketitle	
\tableofcontents

\section{Introduction} \label{sec:intro}
Unlike the nonequivariant stable homotopy theory, the $\RO(\Gmr)$-graded homotopy groups of $\mr{C}_2$-equivariant Eilenberg-MacLane spectrum $\mr{H}\ull{\FF}_2$ is not concentrated in degree zero. In fact, the $\RO(\Cmr_2)$-graded homotopy groups $\pi_\star^{\Cmr_2} \Hmr \ull{\FF}_2$ is nonzero in infinitely many degrees. These groups are completely known as a ring due to unpublished work of Stong and published work of Lewis \cite{LewisCP} (see also \cite{HuKriz, DuggerKR}). Explicitly, 
\[ \MM_{\star} := \pi_\star^{\Cmr_2} \Hmr \ull{\FF}_2 \cong \FF_2[a, u] \oplus \FF_2 \left\lbrace \frac{\Theta}{a^{i} u^{j}}: i, j \geq 0  \right\rbrace\]
consists of a polynomial part generated by $a$ and $u$ in degrees $-\sigma$ and $1 - \sigma$, and an infinitely divisible part led by the element $\Theta$ in degree $2 \sigma-2$. This infinitely divisible part is often called the \emph{negative cone}, whereas the polynomial part is referred to as the \emph{positive cone}.

The element $a$ in the positive cone, often referred to as the Euler class, lies in the image of the equivariant Hurewicz map
\[
\begin{tikzcd}
\Hfrak \colon \pi_{\star}^{\Cmr_2}(\SS_{\Cmr_2}) \arrow[r] & \pi_\star^{\Cmr_2} \Hmr \ull{\FF}_2
\end{tikzcd}
\]
induced by the unit map of $\Hmr\ull{\FF}_2$. In fact, $a \in \pi_{-\sigma}^{\Cmr_2}(\SS_{\Cmr_2})$ originates unstably from the inclusion of fixed points
\[
\begin{tikzcd}
S^0 \arrow[r, hook] & S^{\sigma}
\end{tikzcd}
\]
into the one-point compactification of the sign representation $\sigma$. The fact that only the powers of $a$ in the positive cone lie in the Hurewicz image follows from the structure of the genuine $\Cmr_2$-equivariant Adams spectral sequence. Specifically, the positive cone of $\pi_\star^{\Cmr_2} \mr{H}\ull{\FF}_2$ is isomorphic to the $\mathbb{R}$-motivic cohomology of a point. The Steenrod action on this polynomial ring, originally computed by Voevodsky in the motivic setting and translated to the $\Cmr_2$-equivariant context by Hu and Kriz \cite{HuKriz}, ensures that all other filtration zero classes in the positive cone support nontrivial Adams differentials (see \cite[Proposition 3.2]{DI}). 

Unlike the positive cone, the Hurewicz image in the negative cone remained elusive until recently. In a striking new development, Guo et al. \cite{HurHF} completely determined this image for $\Hmr\ull{\FF}_2$, revealing a profound connection to the classical problem of vector fields on spheres. Specifically, an element of the form $\frac{\Theta}{a^k u^n}$ in the negative cone survives in the equivariant Hurewicz image if and only if the classical sphere $S^n$ admits $k$ linearly independent vector fields. 
\begin{figure}[h]
    \begin{minipage}{1\textwidth}
    \centering
    \includegraphics[scale=.45]{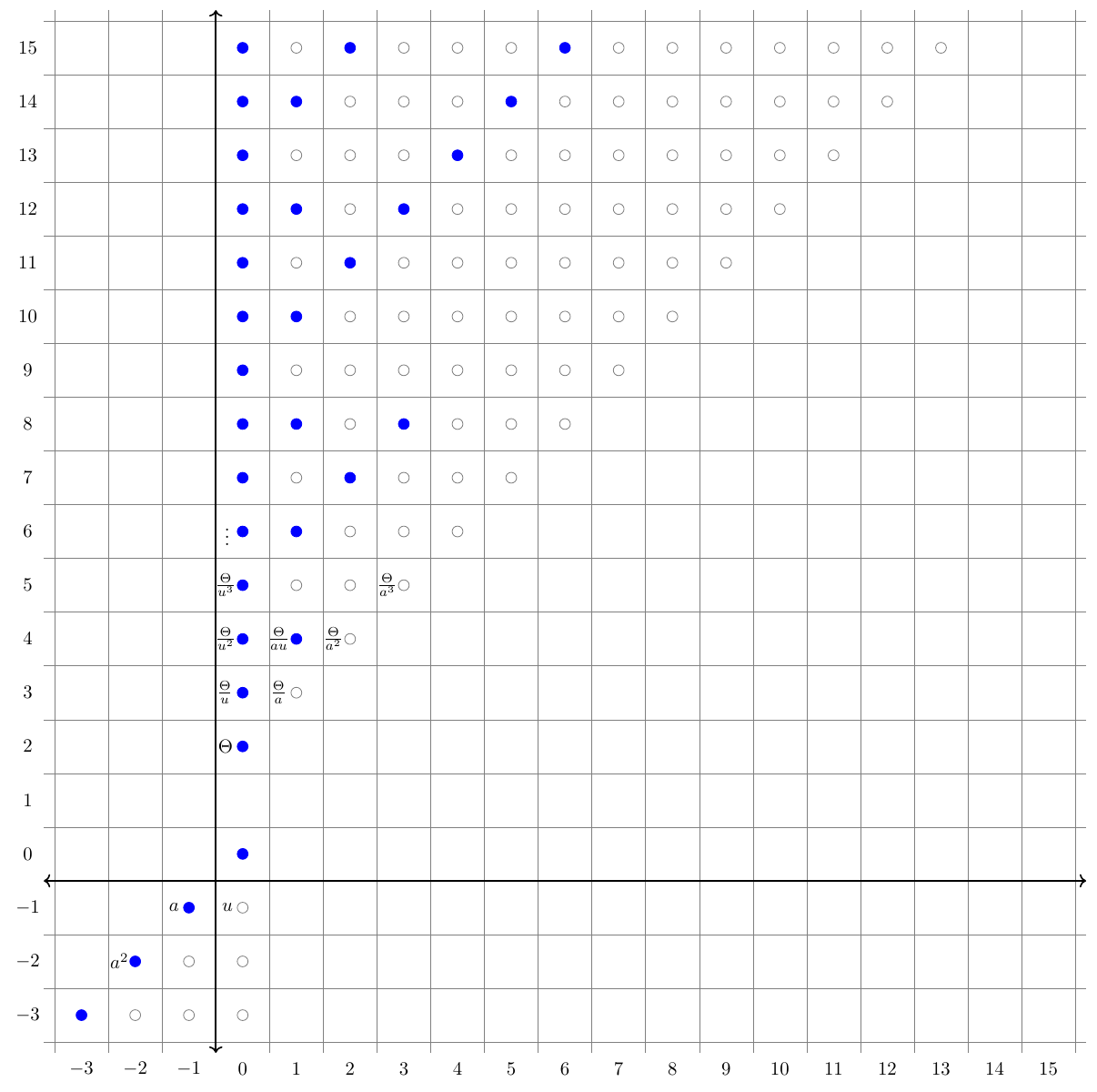}
    \caption{The Hurewicz image of $\MM_2$ in motivic gradings (see \Cref{notn:mot}). }
    \label{fig:HZ}
    \end{minipage}
\end{figure}

\begin{notn} \label{notn:mot}
Throughout the figures in this paper, we pass from the $\RO(\Cmr_2)$-grading $a + b\sigma$ to the motivic bigrading $(s, w)$ (or $(s, t)$), where $s = a + b$ denotes the topological stem and $w = t = b$ denotes the weight.
\end{notn}
\begin{notn}[Vector Fields and James Periodicity] \label{notn:ell_lambda}
Throughout this paper:
\begin{itemize}
\item Let $\ell(n)$ denote the maximal number of linearly independent vector fields on $S^n$. By Adams' theorem \cite{AdamsVF}, writing $n + 1 = (2b+1)2^{c+4d}$ with $b, d \ge 0$ and $0 \le c \le 3$, we have
\[
\ell(n) = 2^c + 8d - 1.
\]

\item Let $\lambda(i)$ denote the $i$-th James number.  Explicitly,
\[
\lambda(i) = 2^{\phi(i-1)} \]
where $\phi(i) = \#\{\, k \in \NN : 1 \le k \le i \text{ and } k \equiv 0, 1, 2, 4 \pmod 8 \,\}$.
\end{itemize}
\end{notn}
\begin{rmk}\label{rmk:ell_lambda}
These two step functions act as order-theoretic inverses of one another, characterized for all $n, i \ge 1$ by the equivalence
\[
\lambda(i) \mid n \quad \Longleftrightarrow \quad i - 1 \le \ell(n-1),
\]
or equivalently, $\ell(n-1) = \max\{\, i - 1 \ge 0 : \lambda(i) \mid n \,\}$ and $\lambda(i) = \min\{\, n \ge 1 : \ell(n-1) \ge i - 1 \,\}$.
\end{rmk}

The primary computational tool in \cite{HurHF} is the two-way correspondence established by Sihao Ma \cite{Ma} between the genuine $\Cmr_2$-equivariant Adams spectral sequence and the Borel $\Cmr_2$-equivariant Adams spectral sequence, and relating them to the classical Adams spectral sequence for stunted real projective spectra. Using this correspondence alongside the Generalized Leibniz Rule and Generalized Mahowald Trick of Li–Wang–Xu \cite{LWX}, they translate the survival of filtration-$0$ classes $\frac{\Theta}{a^k u^n}$ into the classical Adams spectral sequence of stunted projective spaces and deduce that non-surviving filtration-$0$ classes support arbitrarily long genuine $\Cmr_2$-Adams differentials (their Theorem B).

In contrast, our paper uses the Behrens-Shah $u$-self-equivalences \cite{BehShah} of $C(a^k)$ --- which arise as a direct consequence of the classical James periodicity calculations of Adams \cite{AdamsVF} --- to first determine the $\mr{H}\ull{\FF}_2$-Hurewicz image of $C(a^k)$ for all $k$:

\begin{main} \label{main1}
Let $\ell(m)$ denote the maximal number of linearly independent vector fields on $S^m$ and $\lambda(m)$ denote the $m$-th James number (see \Cref{notn:ell_lambda}). For any $k \ge 0$, $0 \le i \le k$, and $n \in \ZZ$, the element
\[
u^{n}a^{k - i} \in \Hmr_{\star}^{\Cmr_2}(C(a^{k+1}); \ull{\FF}_2) \cong \FF_2[a, u^{\pm 1}]/(a^{k+1})
\]
lies in the $\Hmr\ull{\FF}_2$-Hurewicz image of $C(a^{k+1})$ if and only if either of the following equivalent conditions holds:
\begin{enumerate}
    \item \emph{(James periodicity)} The James number $\lambda(i+1)$ divides $n$.
    \item \emph{(Vector fields on spheres)} The sphere $S^{r(n)}$ admits $i$ linearly independent vector fields—that is, $0 \le i \le \ell(r(n))$—where
    \[
    r(n) \in \{0,  \dots, \lambda(k+1)-1 \}
    \]
   is the congruence class of $-n - 1$ modulo $\lambda(k+1)$.
\end{enumerate}
\end{main}

Passing to the limit as $k$ goes to $\infty$, the unit map $\upiota_k \colon \SS_{\Cmr_2} \to C(a^k)$ and the pinch map $\sfp_k \colon C(a^k) \to \Sigma^{1-k\sigma}\SS_{\Cmr_2}$ allow us to completely determine the $\Hmr\ull{\FF}_2$-Hurewicz image of $\SS_{\Cmr_2}$ in both the positive and negative cones of $\MM_{\star}$. This yields an alternative proof of the complete $\Hmr\ull{\FF}_2$-Hurewicz theorem ( see \Cref{main2} below). While our proof share foundational roots with those of \cite{DI, HurHF}, our approach is fundamentally different; for instance, we completely avoid spectral sequences.

\begin{figure}[htbp]
    \begin{minipage}[b]{0.5\textwidth}
    \centering
    \includegraphics[height=10cm, width=1.3\linewidth, keepaspectratio]{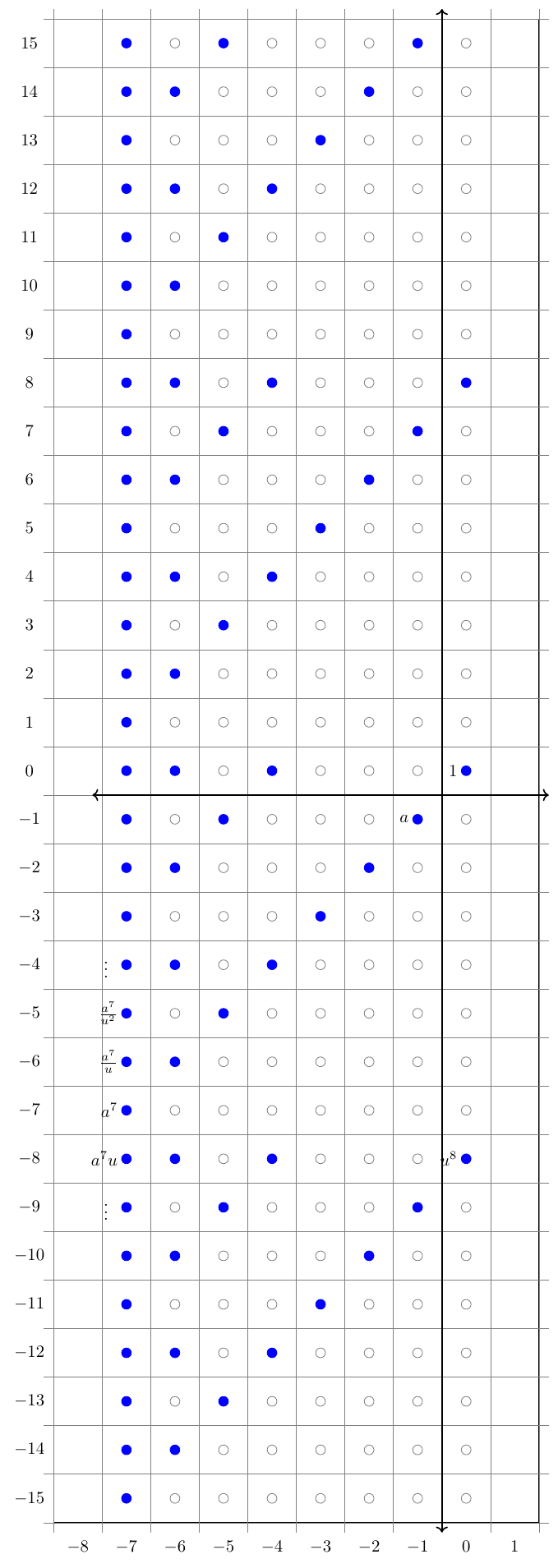}
    \caption{The $\Hmr\ull{\FF}_2$-Hurewicz image of $C(a^{8})$.}
    \end{minipage}%
    \begin{minipage}[b]{0.5\textwidth}
    \centering
    \includegraphics[height=10cm, width=1.3\linewidth, keepaspectratio]{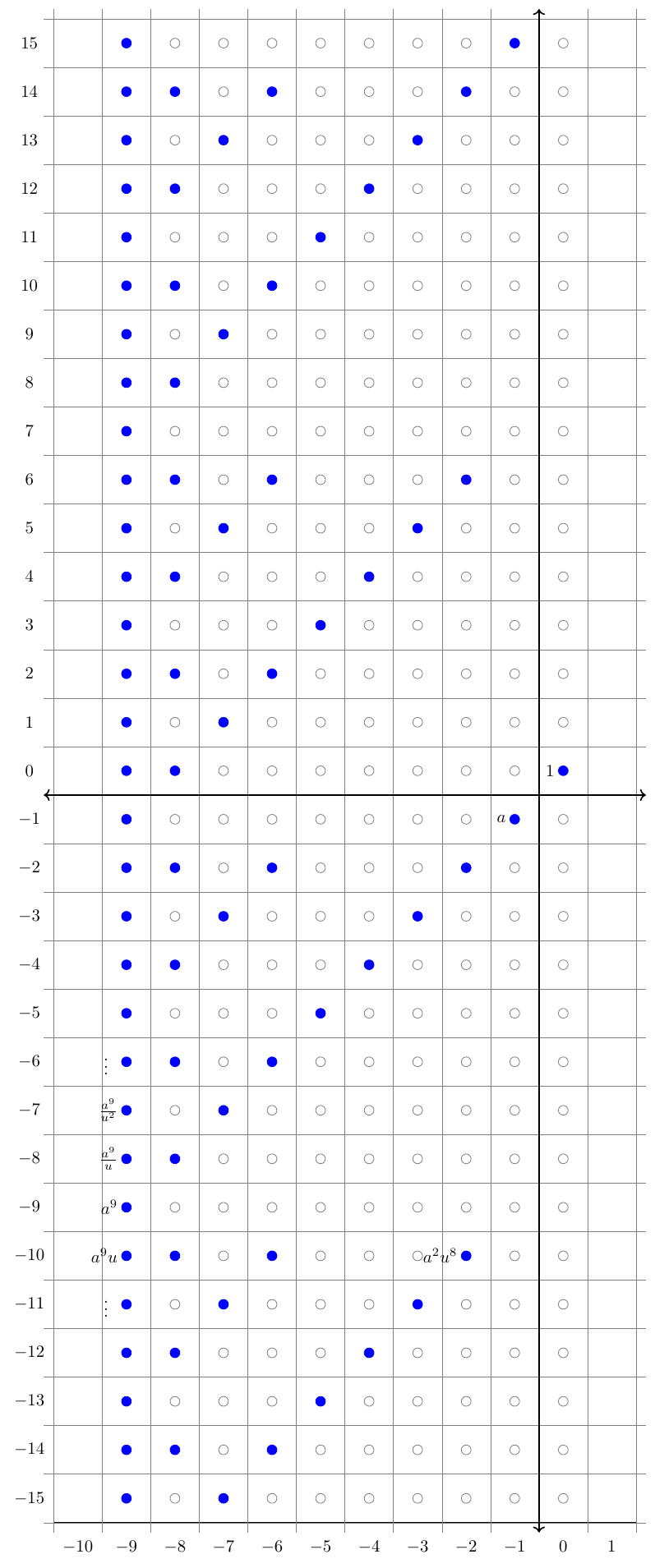}
    \caption{The $\Hmr\ull{\FF}_2$-Hurewicz image of $C(a^{10})$.}
    \end{minipage}
\end{figure}
\begin{main}[{\cite{DI, HurHF}}] \label{main2} 
The $\Hmr\ull{\FF}_2$-Hurewicz image of the $\Cmr_2$-equivariant sphere spectrum $\SS_{\Cmr_2}$ in $\MM_{\star}$ is completely determined as follows:
\begin{enumerate}
    \item \emph{(Positive cone)} For any $k, n \ge 0$, the element 
    \[
    a^k u^n \in \MM_{\star}
    \]
    lies in the $\Hmr\ull{\FF}_2$-Hurewicz image if and only if $n = 0$.
    \item \emph{(Negative cone)} For any $k \ge 0$ and $n \ge 1$, the element 
    \[ 
    \frac{\Theta}{a^k u^{n-1}} \in \MM_{\star} 
    \]
    lies in the $\Hmr\ull{\FF}_2$-Hurewicz image if and only if $S^{n-1}$ admits $k$ linearly independent vector fields (i.e., $0 \le k \le \ell(n-1)$).
\end{enumerate}
\end{main}
\begin{rmk}[Sharp $a$-torsion bounds for Hurewicz lifts] \label{rmk:sharp_a_torsion}
Because any negative-cone class $\frac{\Theta}{a^k u^n} \in \MM_{\star}$ has $a$-torsion order $k+1$ in homology (that is, $a^{k+1} \cdot \frac{\Theta}{a^k u^n} = 0$ while $a^k \cdot \frac{\Theta}{a^k u^n} = \frac{\Theta}{u^n} \neq 0$), any Hurewicz lift $\widetilde{x} \in \pi_{\star}^{\Cmr_2}(\SS_{\Cmr_2})$ must satisfy \[ a^k \cdot \widetilde{x} \neq 0.\] From the perspective of the $\Cmr_2$-equivariant Adams spectral sequence, however, knowing that $\frac{\Theta}{a^k u^n}$ is a permanent cycle only guarantees that $a^{k+1} \cdot \widetilde{x}$ has positive Adams filtration, leaving open the possibility of nontrivial hidden $a$-extensions. A pleasant byproduct of our topological approach is that by constructing Hurewicz lifts directly via the pinch map ${\sf p}_{k+1} \colon C(a^{k+1}) \to \Sigma^{1-(k+1)\sigma}\SS_{\Cmr_2}$, we automatically obtain a lift $\widetilde{x}$ of $\frac{\Theta}{a^k u^n}$ that achieves this sharp $a$-torsion bound at the spectrum level, satisfying $a^{k+1} \cdot \widetilde{x} = 0$ (see \Cref{cor:a_torsion_bound}).
\end{rmk}

\subsection*{Infinite families in $\Cmr_2$-stable stems} \ 

While the Hurewicz image of a ring spectrum directly detects non-trivial elements stable stems, it is also highly effective for identifying elements hidden in the Hurewicz kernel. A notable precedent appears in \cite{BBQ}, where the authors exploited the surjectivity of the $\tmf$-Hurewicz map on $A_1$ \cite{Pham}, together with a $192$-periodic self-map of $A_1$ \cite{BEM}, to construct seven new $192$-periodic elements in chromatic layer $2$ with trivial $\tmf$-Hurewicz image. In a similar vein, our approach utilizes the complete knowledge of the $\mr{H}\ull{\FF}_2$-Hurewicz map on $C(a^k)$, the cofiber of $a^k$, to yield $u$-periodic elements in $\Cmr_2$-equivariant stable stems with trivial $\mr{H}\ull{\FF}_2$-Hurewicz image.
\begin{figure}[h]
    \begin{minipage}{1\textwidth}
    \centering
    \includegraphics[scale=.5]{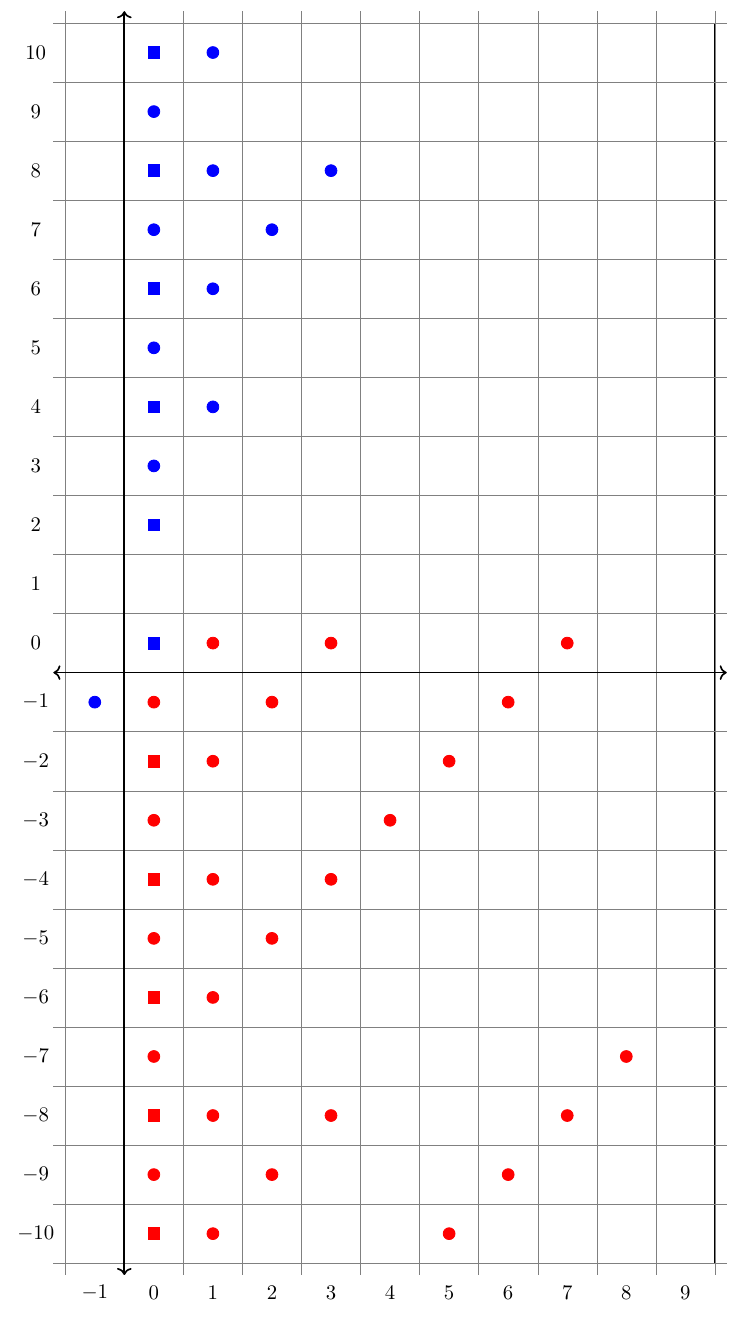}
    \caption{Infinite families in $\pi_\star^{\Cmr_2}(\SS_{\Cmr_2})$ detected by \Cref{main3} in red, and classes detected in $\pi_{\star}^{\Cmr_2}\Hmr\ull{\FF}_2$ in blue. The integral generators are represented by $\blacksquare$ and the $\ZZ/2$-generators are represented by $\bullet$ }
    \label{fig:HZ}
    \end{minipage}
\end{figure}

In this paper, we leverage this complete description of the $\mr{H}\ull{\FF}_2$-Hurewicz image to construct new infinite families in the $\Cmr_2$-equivariant stable stems.
\begin{main} \label{main3} For every $n >0 $, there exists an element $\zeta_{n-1} \in \pi_{(n-1) (1 - \sigma)}^{\Cmr_2}(\SS_{\Cmr_2})$
such that:
\begin{enumerate}
\item  $\zeta_{n-1}$ is $a$-torsion, i.e.,  $ a\cdot \zeta_{n-1} = 0 $
for all $n$. 
\item $\zeta_{n-1}$ is divisible by $a^k$ if and only if $S^{n-1}$ admits $k$ linearly independent vector fields.
\item For any $k \geq 0$ such that $a^k$ divides $\zeta_{n-1}$, any choice of quotient $\frac{\zeta_{n-1}}{a^k}$ has a trivial $\mr{H}\ull{\FF}_2$-Hurewicz image. 
\end{enumerate}
\end{main}
\begin{rmk}
While $\zeta_k$ has trivial $\Hmr\ull{\FF}_2$-Hurewicz image for all $k \in \NN$, the even classes $\zeta_{2k}$ have non-zero image under the $\Hmr\ull{\ZZ}$-Hurewicz homomorphism; in fact, each $\alpha_{2k}$ maps to two times a generator. Furthermore, by \cite{HurKR}, the classes \[ \{ \zeta_{4k+1}, \frac{\zeta_{4k+1}}{a} : k \in \NN \} \] are detected in the Hurewicz image of Atiyah's Real $K$-theory.
\end{rmk}

Notice that the elements $\zeta_{n-1}$ constructed in \Cref{main1} are strictly $a$-torsion. Recall that the $\Cmr_2$-geometric fixed point functor 
\[ 
\begin{tikzcd}
\Phi^{\Cmr_2} \colon \pi_{V}^{\Cmr_2}(\SS_{\Cmr_2}) \rar & \pi_{|V^{\Cmr_2}|}(\SS)
\end{tikzcd}
\]
is implemented by inverting $a$. Because these classes are annihilated by $a$, their localization away from $a$ necessarily vanishes.

On the other hand, consider the underlying restriction map
\[ 
\begin{tikzcd}
\res \colon \pi_{V}^{\Cmr_2}(\SS_{\Cmr_2}) \rar & \pi_{|V|}(\SS),
\end{tikzcd}
\]
which forgets the $\Cmr_2$-action and algebraically corresponds to killing $a$ and setting $u = 1$. Under the Wirthm\"uller isomorphism ($\pi_V^{\Cmr_2}(C(a)) \cong \pi_{|V|}(\SS)$), the map $\res$ fits into the long exact sequence of the cofiber of $a$:
\[
\pi_{V+\sigma}^{\Cmr_2}(\SS_{\Cmr_2}) \xrightarrow{\;\cdot a\;} \pi_V^{\Cmr_2}(\SS_{\Cmr_2}) \xrightarrow{\;\res\;} \pi_{|V|}(\SS).
\]
Thus $\ker(\res) = \mr{Im}(\cdot a)$, meaning that a class in $\pi_V^{\Cmr_2}(\SS_{\Cmr_2})$ restricts non-trivially if and only if it represents a non-zero element in the cokernel of $a$. By \Cref{main3}(2) and \Cref{notn:ell_lambda}, any choice of maximally divided quotient $\frac{\zeta_{n-1}}{a^{\ell(n-1)}}$ is non-trivial in $\mr{coker}(a)$, and therefore restricts to a non-zero element in $\pi_{\ell(n-1)}(\SS)$. We summarize these fixed-point and restriction properties as follows:

\begin{main} \label{main:res_phi}
For every $n > 0$, let $\zeta_{n-1} \in \pi_{(n-1)(1-\sigma)}^{\Cmr_2}(\SS_{\Cmr_2})$ be the element of \Cref{main3}. Then:
\begin{enumerate}
\item For all $0 \le k \le \ell(n-1)$, the $\Cmr_2$-geometric fixed points vanish: 
\[
\Phi^{\Cmr_2}\!\left(\frac{\zeta_{n-1}}{a^k}\right) = 0.
\]
\item For all $0 \le k < \ell(n-1)$, the underlying restriction vanishes: 
\[
\res\!\left(\frac{\zeta_{n-1}}{a^k}\right) = 0.
\]
\item At the maximal divisibility $k = \ell(n-1)$, the underlying restriction of any choice of quotient
\[
\res\!\left(\frac{\zeta_{n-1}}{a^{\ell(n-1)}}\right) \in \pi_{\ell(n-1)}(\SS)
\]
is non-zero. 
\end{enumerate}
\end{main}
\begin{rmk}[The even case $n-1 = 2k$ and $\Hmr\ull{\ZZ}$-detection]
When $n-1 = 2k$ is even, $\ell(2k) = 0$ and $\zeta_{2k}$ itself lies in $\mr{coker}(a)$. In this case, $\zeta_{2k}$ is obtained by pinching a Hurewicz lift of $u^{2k+1} \in \Hmr_{(2k+1)(1-\sigma)}^{\Cmr_2}(C(a); \ull{\FF}_2)$ to $\Sigma^{1-\sigma}\SS_{\Cmr_2}$. Since 
\[
\pi_{(2k+1)(1-\sigma)}^{\Cmr_2}(C(a)) \cong \pi_0(\SS) \cong \ZZ
\]
and the $\Hmr\ull{\FF}_2$-Hurewicz map is reduction modulo $2$, there are $\ZZ$ many choices for this lift. However, passing to $\Hmr\ull{\ZZ}$-homology, the pinched class $\zeta_{2k}$ is detected in the $\Hmr\ull{\ZZ}$-Hurewicz image of $\SS_{\Cmr_2}$ \cite{HurHF}, which fits into the commutative square
\[
\begin{tikzcd}[column sep=large]
\pi_{2k(1-\sigma)}^{\Cmr_2}(\SS_{\Cmr_2}) \rar["\Hfrak_{\ull{\ZZ}}"] \dar["\res"'] & \pi_{2k(1-\sigma)}^{\Cmr_2}(\Hmr\ull{\ZZ}) \cong \ZZ\{u_{2\sigma}^k\} \dar["\res", "\cong"'] \\
\pi_0(\SS) \rar["\cong"', "\Hfrak_{\ZZ}"] & \Hmr_0(\SS; \ZZ) \cong \ZZ.
\end{tikzcd}
\]
Choosing the lift corresponding to a generator of $\Hmr_{(2k+1)(1-\sigma)}^{\Cmr_2}(C(a); \ull{\ZZ}) \cong \ZZ$ sends $\zeta_{2k}$ to $\pm 2 u_{2\sigma}^k \in \pi_{2k(1-\sigma)}^{\Cmr_2}(\Hmr\ull{\ZZ})$, which restricts to $\pm 2 \in \pi_0(\SS) \cong \ZZ$.
\end{rmk}

\begin{rmk}[Indeterminacy and coherent choices]
Neither $\zeta_{n-1}$ nor the quotient $\frac{\zeta_{n-1}}{a^{\ell(n-1)}}$ is unique: Hurewicz lifts are only well-defined modulo $\ker(\Hfrak)$, and $a$-quotients are only defined modulo $\ker\!\left(a^{\ell(n-1)}\right)$. In \Cref{rmk:coherent_lifts}, we compatibly resolve the $a$-divisibility indeterminacy by first fixing a maximally $a$-divisible class $\overline{\zeta}_{n-1}$ and setting $\zeta_{n-1} := a^{\ell(n-1)}\overline{\zeta}_{n-1}$. Even so, the underlying restriction $\res\!\left(\frac{\zeta_{n-1}}{a^{\ell(n-1)}}\right)$ is in general only determined modulo $\res\!\left(\ker(a^{\ell(n-1)})\right)$.
\end{rmk}

\begin{rmk}[Restriction and the image of $J$]
When $n-1$ is odd, $\ell := \ell(n-1) \ge 1$, so $\res(\zeta_{n-1}) = 0$ while $\res\!\left(\frac{\zeta_{n-1}}{a^{\ell}}\right) \neq 0 \in \pi_{\ell}(\SS)$. For $16 \nmid n$ (where $\ell = 2^c - 1 \in \{1, 3, 7\}$), $\zeta_{n-1}$ is detected by $a^{2^c-1}\tau^{n-1}h_c$ in \cite{GI}, so $\frac{\zeta_{n-1}}{a^{2^c-1}}$ is detected by $\tau^{n-1}h_c$ and restricts to an odd multiple of $\eta, \nu,$ or $\sigma$. More generally, by choosing the lift $\widetilde{\upzeta}_{0,n}^{(\ell+1)}$ via the Behrens--Shah equivalence \cite{BehShah} and Atiyah's Thom spectrum formulation of $J$ \cite{AtiyahThom}, one can always arrange for $\res(\overline{\zeta}_{n-1})$ to generate the $2$-primary image of $J$ in $\pi_{\ell}(\SS)$, though an arbitrary quotient $\frac{\zeta_{n-1}}{a^{\ell}}$ may differ by $\res(\ker(a^{\ell}))$ (which for $\ell \ge 8$ can contain classes in $\mr{coker}(J)$).
\end{rmk}

\subsection*{Acknowledgments}
The authors are grateful to Mark Behrens for an illuminating conversation at the International Workshop on Algebraic Topology (IWoAT). While we had already envisioned the construction of the infinite families in the Hurewicz kernel, his intuition that the results of \cite{HurHF} could be tied to the $u$-self-equivalences of \cite{BehShah} proved crucial to our determination of the $\Hmr\ull{\FF}_2$-Hurewicz image of $C(a^k)$. 

We also thank Sihao Ma, Zhouli Xu, and Zhangjie Zhang for valuable discussions, and the organizers of IWoAT for providing a stimulating research environment.

\subsection*{Organization of the paper.}
In \Cref{sec:Hurewicz_Cak}, we review  the Behrens--Shah $u$-self-equivalences and determine the $\Hmr\ull{\FF}_2$-Hurewicz image of $C(a^k)$ for all $k$, our \Cref{main1}. 

  In \Cref{sec:Hurewicz_S}, we use \Cref{main1} to determine the $\Hmr\ull{\FF}_2$-Hurewicz image of $\SS_{\Cmr_2}$, proving \Cref{main2}. 
  
  Finally, in \Cref{sec:infinite_families},  we identify new infinite families in $\Cmr_2$-stable stems and complete the proof of \Cref{main3}.

\section{The $\Hmr\ull{\FF}_2$-Hurewicz image of $\Cmr(a^k)$} \label{sec:Hurewicz_Cak}

At the level of pointed $\Cmr_2$-spaces, the unit sphere $S(k\sigma)$ and unit disk $D(k\sigma)$ of the representation $k\sigma$ fit into a cofiber sequence 
\begin{equation} \label{eq:unstable_cofiber}
\begin{tikzcd} 
S(k \sigma)_+ \rar["c"] & S^{0} \rar[hook, "a_s^{k}"] & S^{k \sigma},
\end{tikzcd}
\end{equation}
where $c$ collapses $S(k\sigma)$ to the non-basepoint of $S^0 \simeq D(k\sigma)_+$, and $a_s^k \colon S^0 \hookrightarrow S^{k \sigma}$ denotes the unstable inclusion of the $\Cmr_2$-fixed points. We use the notation $a_s^k$ to distinguish this unstable map from the stable Euler class introduced in \Cref{sec:intro}:
\[ 
\begin{tikzcd}
a^{k} := \Sigma^{\infty - k \sigma}_{\Cmr_2} a_s^k \colon \Sigma^{-k \sigma}\SS_{\Cmr_2} \rar & \SS_{\Cmr_2}. 
\end{tikzcd}
\]
Equivalently, it is fruitful to view $a^k$ as the Spanier--Whitehead dual of $a_s^k$. Applying the Spanier--Whitehead duality functor \[ \SS_{\Cmr_2}^{(-)} := F(\Sigma^\infty_{\Cmr_2}(-), \SS_{\Cmr_2})\] to \eqref{eq:unstable_cofiber} yields the defining cofiber sequence of $C(a^k)$ in $\Sp^{\Cmr_2}$:
\begin{equation} \label{eq:stable_cofiber}
\begin{tikzcd}
\Sigma^{-k\sigma}\SS_{\Cmr_2} \rar["a^k"] & \SS_{\Cmr_2} \rar["\iota"] & C(a^k),
\end{tikzcd}
\end{equation}
thereby identifying $C(a^k)$ with the Spanier--Whitehead dual of $S(k\sigma)_+$:
\begin{equation} \label{eq:Cakasfunction}
C(a^k) \simeq \SS_{\Cmr_2}^{S(k\sigma)_+}. 
\end{equation}
This identification immediately endows $C(a^k)$ with a highly structured multiplication.
\begin{lem} \label{lem:Cak_Einfty}
For every $k \ge 1$, the cofiber spectrum $C(a^k) \simeq \SS_{\Cmr_2}^{S(k\sigma)_+}$ is a genuine $\Cmr_2$-equivariant $\EE_\infty$ ring spectrum. 
\end{lem}

\begin{proof}
As a  pointed $\Cmr_2$-spaces, $S(k\sigma)_+$ is a cocommutative coalgebra with diagonal \[ \Delta_+ \colon S(k\sigma)_+ \to S(k\sigma)_+ \sma S(k\sigma)_+\]  and counit \[ c \colon S(k\sigma)_+ \to S^0.\] 
Since the contravariant function spectrum functor 
\[ 
\begin{tikzcd}
\SS_{\Cmr_2}^{(-)} \colon (\Top_*^{\Cmr_2})^{\mr{op}} \rar & \Sp^{\Cmr_2}
\end{tikzcd}\]  is lax $\Cmr_2$-symmetric monoidal, it sends cocommutative coalgebras to $\Cmr_2$-commutative algebras, equipping $C(a^k) \simeq \SS_{\Cmr_2}^{S(k\sigma)_+}$ with the structure of a genuine $\Cmr_2$-equivariant $\EE_\infty$ ring spectrum. 
\end{proof}

By \Cref{lem:Cak_Einfty}, the $\RO(\Cmr_2)$-graded Bredon homology of $C(a^k)$ is a ring, which was computed  by Clover May:
\begin{thm}[Clover May {\cite[$\mathsection$3.2]{Clover}}] \label{thm:clover}
For all $k\geq 1$,  
\[ \Hmr_{\star}^{\Cmr_2}(C(a^k); \ull{\FF}_2) \cong \Hmr^{-\star}_{\Cmr_2}(S(k\sigma)_+; \ull{\FF}_2) \cong \FF_2[a, u^{\pm 1}]/(a^{k}).
\]
Here, the generators $a$ and $u$ denote the images of the corresponding classes in $\MM_{\star}$ under the unit map $\iota$ of \eqref{eq:stable_cofiber}.
\end{thm}

\subsection{The Behrens-Shah $u$-self-equivalences} \ 

A notable feature of the homology rings $\Hmr_{\star}^{\Cmr_2}(C(a^k); \ull{\FF}_2)$ is that $u$ is invertible for all $k \geq 1$, in stark contrast to $\MM_{\star}$. Topologically, this algebraic invertibility is explained by the existence of periodic $u$-self-equivalences on $C(a^k)$. 

When $k = 1$, we have $S(\sigma) = \Cmr_2$, so $C(a) \simeq \SS_{\Cmr_2}^{C_{2+}}$. In this case, the class \[u \in \Hmr_{1-\sigma}^{\Cmr_2}(C(a); \ull{\FF}_2)\] itself lifts to a self-equivalence of $C(a)$. To see this geometrically, observe that a shearing map induces a $\Cmr_2$-equivariant homeomorphism
\[
\Sigma^1 C_{2+} \cong \Sigma^\sigma C_{2+}.
\]
Applying the Spanier--Whitehead duality functor $\SS_{\Cmr_2}^{(-)}$  yields a $\Cmr_2$-equivariant self-equivalence
\[
\begin{tikzcd}
\upnu_1 \colon  C(a) \rar["\simeq"] & \Sigma^{1-\sigma}C(a)
\end{tikzcd}
\]
that acts on $\Hmr_{\star}^{\Cmr_2}(C(a); \ull{\FF}_2) \cong \FF_2[u^{\pm 1}]$ as multiplication by $u$.

For $k > 1$, the class $u$ no longer lifts to a self-map of $C(a^k)$, but Behrens and Shah showed that a suitable power of $u$ does:
\begin{thm}[Behrens--Shah {\cite[Theorem 7.7]{BehShah}}] \label{thm:BS} 
For every $k \ge 1$, there exists a $\Cmr_2$-equivariant self-equivalence 
\[
\begin{tikzcd}
\upnu_k \colon  C(a^k) \rar["\simeq"] & \Sigma^{\lambda(k)(1- \sigma )}C(a^k)
\end{tikzcd}
\] 
that induces multiplication by $u^{-\lambda(k)}$ on $\Hmr_{\star}^{\Cmr_2}(C(a^k); \ull{\FF}_2)$, where $\lambda(k)$ denotes the $k$-th James number defined in \Cref{notn:ell_lambda}. 
\end{thm}

Let $\Sp_{(2)}$ and $\Sp^{\Cmr_2}_{(2)}$ denote the categories of $2$-local spectra and $2$-local $\Cmr_2$-spectra, respectively. In this paper, we require a sharp converse to \Cref{thm:BS}:

\begin{thm} \label{thm:iffJames} 
For any $k \ge 1$ and $r \in \ZZ$, the homology isomorphism induced by multiplication by $u^r$,
\[
\begin{tikzcd}[column sep=large]
\cdot u^{r} \colon \Hmr_{\star - r(1- \sigma)}^{\Cmr_2}(C(a^k); \ull{\FF}_2) \rar["\cong"] & \Hmr_{\star }^{\Cmr_2}(C(a^k); \ull{\FF}_2),
\end{tikzcd}
\]
lifts to a self-equivalence $\upnu \colon \Sigma^{r(1- \sigma)} C(a^k) \xrightarrow{\simeq} C(a^k)$ in $\Sp^{\Cmr_2}$ (or $2$-locally in $\Sp^{\Cmr_2}_{(2)}$) if and only if $\lambda(k) \mid r$. 
\end{thm}

To prove this stronger result, let us briefly review the geometric input behind \Cref{thm:BS}. Because $C(a^k)$ is the Spanier--Whitehead dual of the free $\Cmr_2$-space $S(k \sigma)_+$, constructing a $\Cmr_2$-equivariant self-equivalence 
\begin{equation} \label{eq:eq_r}
 \Sigma^{r} C(a^k) \simeq  \Sigma^{r\sigma}C(a^k)
 \end{equation}
is equivalent to a non-equivariant stable equivalence at the level of $\Cmr_2$-fixed points:
\begin{equation} \label{eqn:noneqequiv}
\begin{tikzcd}
F(\Sigma^r\RR\PP^{k-1}_+, \SS) \ar[simeq, d] \ar[rr, "\simeq"] && F(P_{r}^{r+k-1}, \SS)  \\
 F\left( (\Sigma^r S(k\sigma)_+)/_{\Cmr_2}, \SS_{\Cmr_2} \right) \dar[simeq ] && F\left( (\Sigma^{r\sigma} S(k\sigma)_+)/_{\Cmr_2}, \SS_{\Cmr_2}\right) \uar[simeq]\\
  F(\Sigma^r S(k\sigma)_+, \SS_{\Cmr_2})^{\Cmr_2} \ar[simeq, rr] && F(\Sigma^{r\sigma} S(k\sigma)_+, \SS_{\Cmr_2})^{\Cmr_2} \uar[simeq]
\end{tikzcd}
\end{equation}
where $P_r^{r+k-1} = \RR\PP^{r+k-1}/\RR\PP^{r-1}$ denotes the stunted real projective space. 

The question of determining the exact values of $r$ for which $P_r^{r+k-1}$ is stably equivalent to $\Sigma^r \RR\PP^{k-1}_+$ was raised by James \cite{James} and Atiyah \cite{AtiyahThom}, and completely resolved by Adams \cite{AdamsVF} in his solution to the vector fields on spheres problem:

\begin{thm}[ J.F. Adams {\cite[Theorem 1.2 and Theorem 7.4]{AdamsVF}}] \label{thm:Adams_periodicity}
For any $k \ge 1$ and $r \ge 1$, the following are equivalent:
\begin{enumerate}
    \item The stunted real projective space $P_r^{r+k-1}$ is stably homotopy equivalent to $\Sigma^r \RR\PP^{k-1}_+$ (integrally or $2$-locally).
    \item The bottom cell $S^r \hookrightarrow P_r^{r+k-1}$ splits off stably (or dually, the top cell of $D(P_r^{r+k-1})$ splits off stably).
    \item The integer $r$ is a multiple of the James number $\lambda(k) = 2^{\phi(k-1)}$.
\end{enumerate}
\end{thm}
\begin{rmk} \label{rmk:Adams_translation}
While the implications $(3) \implies (1) \implies (2)$ in \Cref{thm:Adams_periodicity} follow integrally from Atiyah's Thom isomorphism $P_r^{r+k-1} \cong \Th(r\LL^{(k-1)})$ \cite{AtiyahThom} (where $\LL^{(k-1)}$ is the tautological line bundle over $\RR\PP^{k-1}$) and Adams' calculation \cite[Theorem 7.4]{AdamsVF}
\[ \widetilde{\KO}(\RR\PP^{k-1}) \cong \ZZ/\lambda(k), \] 
the converse $(2) \implies (3)$ (even $2$-locally) is a direct translation of \cite[Theorem 1.2]{AdamsVF} via \Cref{rmk:ell_lambda}. Specifically, in Adams' notation $\rho(r) = \ell(r-1) + 1$, \cite[Theorem 1.2]{AdamsVF} uses the $2$-group $\widetilde{\KO}(P_r^{r+\rho(r)})$ to show that the bottom cell $S^r \hookrightarrow P_r^{r + \rho(r)}$ does not split off even $2$-locally (i.e., there is no map $P_r^{r+\rho(r)} \to S^r$ of odd degree on the bottom cell). If $\lambda(k) \nmid r$, then \Cref{rmk:ell_lambda} gives 
\[k - 1 \ge \ell(r-1) + 1 = \rho(r).\] 
So $P_r^{r+\rho(r)}$ is a subcomplex of $P_r^{r+k-1}$ containing the bottom cell $S^r$. Consequently, any $2$-local stable retraction of $P_r^{r+k-1}$ onto $S^r$ would restrict to a $2$-local retraction of $P_r^{r+\rho(r)}$ onto $S^r$, which is impossible.
\end{rmk}

\bigskip
\begin{proof}[{\bf Proof of \Cref{thm:BS}}] 
It follows from \Cref{eqn:noneqequiv} and \Cref{thm:Adams_periodicity}, that there exists a $\Cmr_2$-equivariant stable equivalence as in \eqref{eq:eq_r} leading to the construction of $\upnu_k$. 

Next, we want to show that the $\upnu_k$ induces multiplication by  $u^{-\lambda(k)}$ in $\Hmr\ull{\FF}_2$-homology. Since $C(a^k)$ is a ring spectrum it is enough to show that the homotopy class of the composite map 
\begin{equation} \label{eq:split_summand_composite}
\begin{tikzcd}
\Sigma^{-r}\SS_{\Cmr_2} \rar["\Sigma^{-r}\iota"] & \Sigma^{-r} C(a^k) \rar["\simeq"', "\Sigma^{r \sigma}\upnu_k"] & \Sigma^{-r\sigma} C(a^k) \rar["1 \sma \eta"] & \Sigma^{-r\sigma} C(a^k) \sma \Hmr\ull{\FF}_2
\end{tikzcd}
\end{equation}
is $u^{-r}$, where $r = \lambda(k)$.

Since $\Sigma^{-r\sigma} C(a^k) \simeq F(\Sigma^{r\sigma} S(k\sigma)_+, \SS_{\Cmr_2})$ and $\Cmr_2$ acts freely on $\Sigma^{r\sigma} S(k\sigma)_+$ with quotient $P_r^{r+k-1}$, we have canonical isomorphisms
\begin{eqnarray*}
\left[\Sigma^{-r}\SS_{\Cmr_2},\, \Sigma^{-r\sigma} C(a^k) \sma \Hmr\ull{\FF}_2\right]^{\Cmr_2} &\cong& \Hmr^r_{\Cmr_2}\!\left(\Sigma^{r\sigma} S(k\sigma)_+; \ull{\FF}_2\right) \\
&\cong& \Hmr^r\!\left(P_r^{r+k-1}; \FF_2\right) \\ 
&\cong& \Hmr_{-r}\!\left(D(P_r^{r+k-1}); \FF_2\right).
\end{eqnarray*}

On $\Cmr_2$-fixed points, the unit map $\Sigma^{-r}\iota$ is identified with $D(\Sigma^r c) \colon S^{-r} \to D(\Sigma^r \RR\PP^{k-1}_+)$, which is the inclusion of a split summand because the collapse map 
\[ \Sigma^r c \colon \Sigma^r \RR\PP^{k-1}_+  \to S^r\] 
splits the bottom cell. Consequently, the composite of the first two maps in \eqref{eq:split_summand_composite} induces on $\Cmr_2$-fixed points the inclusion of a split summand
\[
\begin{tikzcd}
S^{-r} \rar["D(\Sigma^r c)"] & D(\Sigma^r \RR\PP^{k-1}_+) \rar["\simeq"] & D(P_r^{r+k-1}),
\end{tikzcd}
\]
which is dual to the stable retraction $P_r^{r+k-1} \to S^r$ splitting off the bottom cell. Because a split inclusion induces an injection on $\Hmr^{r}(-; \FF_2)$, the composite \eqref{eq:split_summand_composite} represents the unique nonzero class
\[
u^{-r} \in \Hmr_{-r(1-\sigma)}^{\Cmr_2}(C(a^k); \ull{\FF}_2) \cong \Hmr^r(P_r^{r+k-1}; \FF_2) \cong \FF_2
\]
as desired. 
\end{proof} 

\bigskip
\begin{proof}[{\bf Proof of \Cref{thm:iffJames}}] 
If $\lambda(k) \mid r$, write $r = m\lambda(k)$ for $m \in \ZZ$. Then the $m$-fold iterate $\upnu_k^{\circ m}$ of the equivalence from \Cref{thm:BS} (or its inverse) provides an integral self-equivalence $\Sigma^{r(1-\sigma)} C(a^k) \xrightarrow{\simeq} C(a^k)$ in $\Sp^{\Cmr_2}$ inducing multiplication by $u^r$ on $\Hmr_{\star}^{\Cmr_2}(C(a^k); \ull{\FF}_2)$.

Conversely, suppose there exists a self-equivalence $\upnu \colon \Sigma^{r(1-\sigma)} C(a^k) \xrightarrow{\simeq} C(a^k)$ in $\Sp^{\Cmr_2}_{(2)}$. Composing $\upnu$ with $\upnu_k^{\circ N}$ for $N \gg 0$ if necessary, we may assume $r \ge 1$. Passing to $\Cmr_2$-fixed points via \eqref{eqn:noneqequiv} and taking Spanier--Whitehead duals yields a $2$-local stable equivalence
\[ 
\Sigma^{r} \RR\PP^{k-1}_+ \simeq_{(2)} P_{r}^{r + k-1},
\]
which by \Cref{thm:Adams_periodicity} implies $\lambda(k) \mid r$.  
\end{proof}
\subsection{Proof of \Cref{main1} } \ 

The following is an immediate corollary of \Cref{thm:iffJames}:
\begin{cor} \label{cor:u_Hurewicz}
For any $k \ge 1$ and $r \in \ZZ$, the element \[ u^{r} \in \Hmr_{r(1-\sigma)}^{\Cmr_2}(C(a^k); \ull{\FF}_2)\] lies in the $\Hmr\ull{\FF}_2$-Hurewicz image of $C(a^k)$ if and only if the $k$-th James number $\lambda(k)$ divides $r$.  
\end{cor}

\begin{proof}
If $\lambda(k) \mid r$, then by \Cref{thm:iffJames} there is a self-equivalence  in $\Sp^{\Cmr_2}$
\[ 
\begin{tikzcd}
\upnu : \Sigma^{r(1-\sigma)} C(a^k) \rar["\simeq"] & C(a^k)
\end{tikzcd} 
\] 
 inducing multiplication by $u^r$ on $\Hmr\ull{\FF}_2$-homology. Precomposing $\upnu$ with the suspended unit map $\Sigma^{r(1-\sigma)}\iota \colon \Sigma^{r(1-\sigma)}\SS_{\Cmr_2} \to \Sigma^{r(1-\sigma)} C(a^k)$ gives a Hurewicz lift of $u^r$.

Conversely, suppose $u^r$ admits a Hurewicz lift $\tilde{u} \colon \Sigma^{r(1-\sigma)}\SS_{\Cmr_2} \to C(a^k)$. Using the ring multiplication $\mu \colon C(a^k) \sma C(a^k) \to C(a^k)$ from \Cref{lem:Cak_Einfty}, we form the composite self-map
\[
\begin{tikzcd}
\upnu \colon \Sigma^{r(1-\sigma)} C(a^k) \rar["\tilde{u} \sma 1"] & C(a^k) \sma C(a^k) \rar["\mu"] & C(a^k).
\end{tikzcd}
\]
Because the $\Hmr\ull{\FF}_2$-Hurewicz map is a ring homomorphism, $\upnu$ induces multiplication by $u^r$ on $\Hmr_{\star}^{\Cmr_2}(C(a^k); \ull{\FF}_2)$, which is an isomorphism. Since $C(a^k)$ is a finite free $\Cmr_2$-spectrum, an $\Hmr\ull{\FF}_2$-homology isomorphism is a $2$-local equivalence in $\Sp^{\Cmr_2}_{(2)}$
. By \Cref{thm:iffJames}, it follows that $\lambda(k) \mid r$.
\end{proof}

\bigskip
\begin{proof}[{\bf Proof of \Cref{main1}}] We first show by induction on $k \ge 0$ that \[ u^n a^{k-i} \in \Hmr_{\star}^{\Cmr_2}(C(a^{k+1}); \ull{\FF}_2)\] lies in the $\Hmr\ull{\FF}_2$-Hurewicz image if and only if $\lambda(i+1) \mid n$.

The base case $k = 0$ (where $i = 0$) is immediate from \Cref{cor:u_Hurewicz}, since \[ \lambda(1) = 2^{\phi(0)} = 1\] divides every $n \in \ZZ$. For the inductive step, assume the claim holds for $C(a^k)$ with $k \ge 1$. When $i = k$, \Cref{cor:u_Hurewicz} directly implies that $u^n \in \Hmr_{\star}^{\Cmr_2}(C(a^{k+1}); \ull{\FF}_2)$ lies in the Hurewicz image if and only if $\lambda(k+1) \mid n$. 

It remains to treat the range $0 \le i < k$. Consider the map \[ \alpha^{(k)} :\Sigma^{-\sigma} C(a^k) \to C(a^{k+1}) \] induced by the commutative diagram of horizontal cofiber sequences:
\begin{equation} \label{eq:alpha}
\begin{tikzcd}
\Sigma^{-(k+1)\sigma}\SS_{\Cmr_2} \rar["a^{k}"] \dar[equal] & \Sigma^{- \sigma}\SS_{\Cmr_2}\rar["\Sigma^{-\sigma}\iota"] \dar["a"] & \Sigma^{-\sigma}C(a^{k}) \dar[dashed, "\alpha^{(k)}"] \\
\Sigma^{-(k+1)\sigma}\SS_{\Cmr_2} \rar["a^{k+1}"'] & \SS_{\Cmr_2} \rar["\iota"] & C(a^{k+1}). 
\end{tikzcd}
\end{equation}
 In $\Hmr\ull{\FF}_2$-homology, $\alpha^{(k)}$ induces multiplication by $a$, i.e., the homomorphism
\[ 
\begin{tikzcd}
\alpha_{\ast}^{(k)} \colon \Hmr_{\star + \sigma}^{\Cmr_2}(C(a^{k}); \ull{\FF}_2) \dar[cong] \ar[rrr] &&& \Hmr_{\star}^{\Cmr_2}(C(a^{k+1}); \ull{\FF}_2) \dar[cong] \\
\FF_2[a, u^{\pm 1}]/(a^k)  &&& \FF_2[a, u^{\pm 1}]/(a^{k+1})
\end{tikzcd}
\]
maps $u^{n}a^{k-1-i} \mapsto u^{n}a^{k-i}$.

Next, we study the commutative diagram induced by $\alpha^{(k)}$ and the naturality of the Hurewicz map:
\[ 
\begin{tikzcd}
\pi_{V + \sigma}^{\Cmr_2}(C(a^k)) \rar[->>, "\alpha_{\ast}^{(k)}"] \dar["\Hfrak"'] & \pi_{V}^{\Cmr_2}(C(a^{k+1})) \dar["\Hfrak"] \\
\Hmr_{V + \sigma}^{\Cmr_2}(C(a^{k}); \ull{\FF}_2) \rar["\cong"', "\alpha_{\ast}^{(k)}"] & \Hmr_{V}^{\Cmr_2}(C(a^{k+1}); \ull{\FF}_2)
\end{tikzcd}
\]
We claim that the top horizontal map in the diagram above is a surjection whenever the underlying dimension $|V|$ of $V = \deg(u^n a^{k-i})$ is negative. Indeed, by the octahedral axiom, the cofiber of $\alpha^{(k)}$ is $C(a) \simeq \SS_{\Cmr_2}^{\Cmr_{2+}}$, yielding an exact sequence
\[
\begin{tikzcd}
\pi_{V + \sigma}^{\Cmr_2}(C(a^k)) \rar["\alpha_{\ast}^{(k)}"] & \pi_{V}^{\Cmr_2}(C(a^{k+1})) \rar & \pi_{V}^{\Cmr_2}(C(a)) \cong \pi_{|V|}^{e}(\SS).
\end{tikzcd}
\]
Since the non-equivariant sphere spectrum $\SS$ is connective, $\pi_{|V|}^{e}(\SS) = 0$ whenever $|V| < 0$, so $\alpha_{\ast}^{(k)}$ is surjective on $\pi_{V+\sigma}^{\Cmr_2}$ as claimed.

Because the bottom horizontal map $\alpha_{\ast}^{(k)}$ in the square is an isomorphism sending $u^n a^{k-1-i} \mapsto u^n a^{k-i}$, the surjectivity of the top horizontal map implies that, in this negative underlying dimension range, $u^n a^{k-i}$ lies in the $\Hmr\ull{\FF}_2$-Hurewicz image of $C(a^{k+1})$ if and only if $u^n a^{k-1-i}$ lies in the $\Hmr\ull{\FF}_2$-Hurewicz image of $C(a^k)$. By the inductive hypothesis, this holds if and only if $\lambda(i+1) \mid n$.

Finally, to extend this equivalence to arbitrary $n \in \ZZ$, we apply a suitable iteration, say  $(\upnu_{k+1})^{\circ m}$, of the Behrens--Shah self-equivalence so that the shifted class $u^{n - m\lambda(k+1)} a^{k-i}$ lies in the negative underlying dimension range above. Because $\upnu_{k+1}$ is a $\Cmr_2$-equivariant equivalence inducing multiplication by $u^{-\lambda(k+1)}$ in $\Hmr\ull{\FF}_2$-homology, $u^n a^{k-i}$ lies in the Hurewicz image of $C(a^{k+1})$ if and only if $u^{n - m\lambda(k+1)} a^{k-i}$ does. Since $\lambda(i+1) \mid \lambda(k+1)$ for all $0 \le i < k$, we have
\[
\lambda(i+1) \mid (n - m\lambda(k+1)) \quad \Longleftrightarrow \quad \lambda(i+1) \mid n,
\]
which establishes the result for all $0 \leq i < k$ and $n \in \ZZ$. When $i = k$, \Cref{cor:u_Hurewicz} directly shows that $u^n \in \Hmr_{\star}^{\Cmr_2}(C(a^{k+1}); \ull{\FF}_2)$ lies in the Hurewicz image if and only if $\lambda(k+1) \mid n$, completing the induction. 

To see the equivalence of conditions~(1) and~(2), note that \[ r(n) + 1 \equiv -n \pmod{\lambda(k+1)}.\] Since $\lambda(i+1)$ divides $\lambda(k+1)$ for all $0 \le i \le k$, we have $\lambda(i+1) \mid n$ if and only if $\lambda(i+1) \mid (r(n) + 1)$, which by \Cref{rmk:ell_lambda} holds if and only if $0 \le i \le \ell(r(n))$.
\end{proof}

\section{The $\Hmr\ull{\FF}_2$-Hurewicz image of $\SS_{\Cmr_2}$} \label{sec:Hurewicz_S}

In this section, we show that \Cref{main1} completely determines the image of the equivariant Hurewicz map 
\[
\begin{tikzcd}
\Hfrak \colon \pi_{\star}^{\Cmr_2}(\SS_{\Cmr_2}) \arrow[r] & \pi_\star^{\Cmr_2} \Hmr \ull{\FF}_2 \cong \MM_\star. 
\end{tikzcd}
\]
We begin by fixing some notation. 

\begin{notn} \label{notn:cofiber_maps}
For the remainder of the paper, we use the following conventions for each $k \ge 1$:
\begin{itemize}
\item Let $\upiota_k$ and ${\sf p}_k$ denote the inclusion of the bottom cell and the pinch map of $C(a^k)$, respectively, which fit into the defining cofiber sequence
\begin{equation} \label{eq:full_cofiber_k}
\begin{tikzcd}
\Sigma^{-k\sigma}\SS_{\Cmr_2} \rar["a^k"] & \SS_{\Cmr_2} \rar["\upiota_k"]  & C(a^k) \rar["{\sf p}_k"] & \Sigma^{1 - k \sigma} \SS_{\Cmr_2}.
\end{tikzcd}
\end{equation}
\item Let $\alpha^{(k)} \colon \Sigma^{- \sigma}C(a^k) \to C(a^{k+1})$ denote the map defined in \eqref{eq:alpha}.
\item Let $\beta^{(k)} \colon C(a^{k+1}) \to C(a^{k})$ denote the canonical reduction map fitting into the commutative diagram of horizontal cofiber sequences:
\begin{equation} \label{eq:beta}
\begin{tikzcd}
\Sigma^{-(k+1)\sigma}\SS_{\Cmr_2} \rar["a^{k+1}"] \dar["a"'] & \SS_{\Cmr_2}\rar["\upiota_{k+1}"] \dar[equal] & C(a^{k+1}) \dar[dashed, "\beta^{(k)}"] \\
\Sigma^{-k\sigma}\SS_{\Cmr_2} \rar["a^{k}"'] & \SS_{\Cmr_2} \rar["\upiota_k"'] & C(a^{k}). 
\end{tikzcd}
\end{equation}
\end{itemize}
\end{notn} 

Applying $\Hmr_{\star}^{\Cmr_2}(-; \ull{\FF}_2)$ to \eqref{eq:full_cofiber_k} yields the short exact sequence
\[
\begin{tikzcd}
0 \rar & \mr{coker}(\cdot a^k) \rar["\upiota_{k\ast}"] & \Hmr_{\star}^{\Cmr_2}(C(a^k); \ull{\FF}_2) \rar["{\sf p}_{k\ast}"] & \ker(\cdot a^k) \rar & 0.
\end{tikzcd}
\]
On $\MM_{\star} = \pi_{\star}^{\Cmr_2}\Hmr\ull{\FF}_2$, the cokernel of multiplication by $a^k$ is spanned by the monomials $a^i u^n$ in the positive cone with $0 \le i < k$ and $n \ge 0$, whereas the kernel of multiplication by $a^k$ is spanned by the elements $\frac{\Theta}{a^i u^n}$ in the negative cone with $0 \le i < k$ and $n \ge 0$. Comparing $\RO(\Cmr_2)$-degrees immediately yields the following description of the induced maps in $\Hmr\ull{\FF}_2$-homology:

\begin{lem} \label{lem:iota_p_homology}
For every $k \ge 1$:
\begin{enumerate}
\item The unit map $\upiota_k \colon \SS_{\Cmr_2} \to C(a^k)$ induces the ring homomorphism
\[ 
\begin{tikzcd}
\upiota_{k \ast} \colon \MM_{\star} \rar & \Hmr_{\star}^{\Cmr_2}(C(a^k); \ull{\FF}_2) \cong \FF_2[a, u^{\pm 1}]/(a^k)
\end{tikzcd}
\]
sends all elements of the negative cone to $0$ and 
\[
\upiota_{k\ast}(a^i u^n) = \begin{cases}
a^i u^n & \text{if } 0 \le i < k \text{ and } n \ge 0, \\
0 & \text{if } i \ge k,
\end{cases}
\]
\item The pinch map ${\sf p}_k \colon C(a^k) \to \Sigma^{1-k\sigma}\SS_{\Cmr_2}$ induces the $\MM_{\star}$-module homomorphism
\[ 
\begin{tikzcd}
{\sf p}_{k \ast} \colon \FF_2[a, u^{\pm 1}]/(a^k) \cong \Hmr_{\star}^{\Cmr_2}(C(a^k); \ull{\FF}_2) \rar & \MM_{\star - 1 + k\sigma}
\end{tikzcd}
\]
given for $0 \le i < k$ and $n \in \ZZ$ by
\[
{\sf p}_{k\ast}(a^{k-1-i} u^n) = \begin{cases}
\dfrac{\Theta}{a^i u^{-n-1}} & \text{if } n < 0, \\[8pt]
0 & \text{if } n \ge 0
\end{cases}
\]

\item The map $\beta^{(k)} \colon C(a^{k+1}) \to C(a^k)$ induces the canonical quotient ring homomorphism
\[
\begin{tikzcd}
\beta^{(k)}_{\ast} \colon \FF_2[a, u^{\pm 1}]/(a^{k+1}) \rar & \FF_2[a, u^{\pm 1}]/(a^k),
\end{tikzcd}
\] 
such that 
\[ a^i u^n \longmapsto \begin{cases}
a^i u^n & \text{if } 0 \le i < k, \\
0 & \text{if } i = k.
\end{cases}
\]
\item The map $\alpha^{(k)} \colon \Sigma^{-\sigma} C(a^k) \to C(a^{k+1})$ of \eqref{eq:alpha} induces the injective $\MM_{\star}$-module homomorphism
\[
\begin{tikzcd}
\alpha^{(k)}_{\ast} \colon \FF_2[a, u^{\pm 1}]/(a^k)  \rar &  \FF_2[a, u^{\pm 1}]/(a^{k+1}) 
\end{tikzcd}
\]
given by multiplication by $a$, i.e. 
\[
a^i u^n \longmapsto a^{i+1} u^n.
\]
for $0 \le i < k$ and $n \in \ZZ$. 
\end{enumerate}
\end{lem}

\bigskip
\begin{proof}[{\bf Proof of \Cref{main2}}] 
We first determine the image of the Hurewicz map on the positive cone. Since $a \in \pi_{-\sigma}^{\Cmr_2}(\SS_{\Cmr_2})$ is a Hurewicz lift of $a \in \MM_{-\sigma}$ and $\Hfrak$ is a ring homomorphism, $a^k = \Hfrak(a^k)$ lies in the Hurewicz image for all $k \ge 0$. 

Conversely, suppose $a^k u^n \in \MM_{\star}$ ($k, n \ge 0$) lies in the $\Hmr\ull{\FF}_2$-Hurewicz image of $\SS_{\Cmr_2}$. By \eqref{eq:beta}, the unit maps $\upiota_m \colon \SS_{\Cmr_2} \to C(a^m)$ and the reduction maps $\beta^{(m)} \colon C(a^{m+1}) \to C(a^m)$ assemble into a commutative tower:
\[ 
\begin{tikzcd}
& \SS_{\Cmr_2} \ar[ddrrr, "\upiota_1"] \ar[ddr, "\upiota_m" description] \ar[dd, "\upiota_{m+1}"']  \\
|[xshift=1.5em]| \dots  \\ 
\dots \rar & C(a^{m+1}) \rar["\beta^{(m)}"'] & C(a^m) \rar & \dots \rar & C(a). 
\end{tikzcd}
\]
For every $m > k$, \Cref{lem:iota_p_homology}(1) gives
\[
\upiota_{m\ast}(a^k u^n) =  a^{(m-1) - (m-1-k)} u^n \neq 0 \in \Hmr_{\star}^{\Cmr_2}(C(a^m); \ull{\FF}_2) \cong \FF_2[a, u^{\pm 1}]/(a^m).
\]
By the naturality of the Hurewicz map along $\upiota_m$, the class $u^n a^{(m-1) - (m-1-k)}$ must lie in the $\Hmr\ull{\FF}_2$-Hurewicz image of $C(a^m)$ for all $m > k$.

Applying \Cref{main1}(1) with $i = m - 1 - k$, we deduce that the James number $\lambda(m - k)$ divides $n$ for every $m > k$. Because $\lim_{m \to \infty} \lambda(m - k) = \infty$, the only non-negative integer $n \ge 0$ divisible by $\lambda(m - k)$ for all $m > k$ is $n = 0$. This proves~(1).

\medskip

Next, we determine the Hurewicz image on the negative cone. Extending the commutative diagram \eqref{eq:alpha} one step to the right shows that the pinch maps ${\sf p}_m$ commute with the maps $\alpha^{(m)}$ forming a dual directed system under $\SS_{\Cmr_2}$:
\[
\begin{tikzcd}
\Sigma^{\sigma - 1} C(a) \rar["\alpha^{(1)}"] \ar[ddrrr, "{\sf p}_1"'] & \dots \rar & \Sigma^{m\sigma - 1} C(a^m) \rar["\alpha^{(m)}"] \ar[ddr, "{\sf p}_m" description] & \Sigma^{(m+1)\sigma - 1} C(a^{m+1}) \rar \ar[dd, "{\sf p}_{m+1}"] & \dots \\
& & & & |[xshift=-1.5em]| \dots \\
& & & \SS_{\Cmr_2}. & 
\end{tikzcd}
\]
Fix $k \ge 0$ and $n \ge 1$, and consider the negative-cone element $\frac{\Theta}{a^k u^{n-1}} \in \MM_V$, which lies in $\RO(\Cmr_2)$-degree
\[
V := (2\sigma - 2) + k\sigma - (n-1)(1 - \sigma) = -(n + 1) + (n + k + 1)\sigma.
\]
If $S^{n-1}$ admits $k$ linearly independent vector fields (so that $0 \le k \le \ell(n-1)$, or equivalently $\lambda(k+1) \mid n$), then for any $m > k$, \Cref{main1} implies that the class
\[
u^{-n} a^{m - 1 - k} \in \Hmr_{V + 1 - m\sigma}^{\Cmr_2}(C(a^m); \ull{\FF}_2)
\]
lies in the $\Hmr\ull{\FF}_2$-Hurewicz image of $C(a^m)$. By \Cref{lem:iota_p_homology}(2), the homological pinch map ${\sf p}_{m\ast}$ sends 
\[u^{-n} a^{m - 1 - k} \mapsto \frac{\Theta}{a^k u^{n-1}}.\] Thus, by the naturality of the Hurewicz map along ${\sf p}_m$, the element $\frac{\Theta}{a^k u^{n-1}}$ lies in the $\Hmr\ull{\FF}_2$-Hurewicz image of $\SS_{\Cmr_2}$.

Conversely, suppose $\frac{\Theta}{a^k u^{n-1}} \in \MM_V$ lies in the $\Hmr\ull{\FF}_2$-Hurewicz image of $\SS_{\Cmr_2}$, lifted by a homotopy class $x \in \pi_V^{\Cmr_2}(\SS_{\Cmr_2})$. We claim that $x$ cannot support an infinite $a$-tower in $\pi_{\star}^{\Cmr_2}(\SS_{\Cmr_2})$---that is, $a^m x = 0$ for $m \gg 0$. Indeed, inverting $a$ computes the homotopy groups of the $\Cmr_2$-geometric fixed points $\Phi^{\Cmr_2}\SS_{\Cmr_2} \simeq \SS$:
\[
\begin{tikzcd}
\Phi^{\Cmr_2} \colon \pi_V^{\Cmr_2}(\SS_{\Cmr_2}) \rar & \pi_{|V^{\Cmr_2}|}(\SS) = \pi_{-(n+1)}(\SS).
\end{tikzcd}
\]
If $x$ generated an infinite $a$-tower, its image $\Phi^{\Cmr_2}(x)$ would be a nonzero element in $\pi_{-(n+1)}(\SS)$. However, because $n \ge 1$, the fixed-point dimension $|V^{\Cmr_2}| = -(n+1) \le -2$ is strictly negative, so $\pi_{-(n+1)}(\SS) = 0$ by the connectivity of $\SS$. 

Consequently, there exists $m > k$ such that $a^m x = 0$. By exactness of the cofiber sequence \eqref{eq:full_cofiber_k}, $x$ lifts along the pinch map ${\sf p}_m$ to a homotopy class $\tilde{x} \in \pi_{V + 1 - m\sigma}^{\Cmr_2}(C(a^m))$ with ${\sf p}_{m\ast}(\tilde{x}) = x$. Consider the commutative naturality square for the Hurewicz map in degree $W := V + 1 - m\sigma = -n(1-\sigma) - (m-1-k)\sigma$:
\[
\begin{tikzcd}[column sep=large]
\pi_{W}^{\Cmr_2}(C(a^m)) \rar["{\sf p}_{m\ast}"] \dar["\Hfrak"'] & \pi_{V}^{\Cmr_2}(\SS_{\Cmr_2}) \dar["\Hfrak"] \\
\Hmr_{W}^{\Cmr_2}(C(a^m); \ull{\FF}_2) \rar["\cong"', "{\sf p}_{m\ast}"] & \MM_{V}.
\end{tikzcd}
\]
Because $-n < 0$, \Cref{lem:iota_p_homology}(2) shows that the bottom horizontal map ${\sf p}_{m\ast}$ is an isomorphism in degree $W$, sending $u^{-n} a^{m-1-k}$ to $\frac{\Theta}{a^k u^{n-1}}$. Commutativity of the square forces $\Hfrak(\tilde{x}) = u^{-n} a^{m-1-k}$, meaning that $u^{-n} a^{m-1-k}$ lies in the $\Hmr\ull{\FF}_2$-Hurewicz image of $C(a^m)$. By \Cref{main1}, this holds only if $\lambda(k+1) \mid (-n)$, or equivalently $0 \le k \le \ell(n-1)$, completing the proof of~(2).
\end{proof}
It is worth emphasizing that our proof of \Cref{main2} via \Cref{main1} also yields a sharp $a$-torsion bound on the Hurewicz lifts in the negative cone.

\begin{cor} \label{cor:a_torsion_bound}
Every element 
\[ 
\frac{\Theta}{a^k u^n} \in \pi_{\star}^{\Cmr_2}\Hmr\ull{\FF}_2 
\]
in the $\Hmr\ull{\FF}_2$-Hurewicz image of $\SS_{\Cmr_2}$ admits a Hurewicz lift $\widetilde{x} \in \pi_{\star}^{\Cmr_2}(\SS_{\Cmr_2})$ such that $a^{k+1} \cdot \widetilde{x} = 0$ (while $a^k \cdot \widetilde{x} \neq 0$). 
\end{cor}

\begin{proof}
By \Cref{main2}(2) and \Cref{main1}, whenever $\frac{\Theta}{a^k u^n}$ lies in the $\Hmr\ull{\FF}_2$-Hurewicz image of $\SS_{\Cmr_2}$, the class 
\[ u^{-(n+1)} \in \Hmr_{\star}^{\Cmr_2}(C(a^{k+1}); \ull{\FF}_2)\] admits a Hurewicz lift $ \widetilde{u} \in \pi_{\star}^{\Cmr_2}(C(a^{k+1})).$ Since
\[ {\sf p}_{k+1\ast}(u^{-(n+1)}) = \frac{\Theta}{a^k u^n}\] by \Cref{lem:iota_p_homology}(2), the image $\widetilde{x} := {\sf p}_{k+1\ast}(\widetilde{u}) \in \pi_{\star}^{\Cmr_2}(\SS_{\Cmr_2})$ is a Hurewicz lift of $\frac{\Theta}{a^k u^n}$. By the exactness of the cofiber sequence \eqref{eq:full_cofiber_k} for $C(a^{k+1})$, the composite $a^{k+1} \circ {\sf p}_{k+1}$ is nullhomotopic, so $a^{k+1} \cdot \widetilde{x} = 0$. On the other hand, \[ a^k \cdot \widetilde{x} \neq 0\] since $\Hfrak(a^k \cdot \widetilde{x}) = \frac{\Theta}{u^n} \neq 0$.
\end{proof}

\begin{rmk}[On the $\Hmr\ull{\ZZ}$-Hurewicz image] \label{rmk:HZ_Hurewicz}
Our method of passing through the Hurewicz image of the cofiber spectra $C(a^k)$ can also be employed to determine the $\Hmr\ull{\ZZ}$-Hurewicz image of $\SS_{\Cmr_2}$. We do not elaborate on the details here, as the proof is essentially identical and the result has already been established by Guo et al.\ \cite{HurHF}.
\end{rmk}

\section{New $\Cmr_2$-equivariant infinite families} \label{sec:infinite_families}

A striking consequence of \Cref{main1} and \Cref{main2} is a fundamental tension between the behavior of positive powers of $u$ in homology and in homotopy. Specifically, for $0 \le i < k$ and $n > 0$, \Cref{main1} shows that the element
\[
a^i u^n \in \Hmr_{\star}^{\Cmr_2}(C(a^k); \ull{\FF}_2) \cong \FF_2[a, u^{\pm 1}]/(a^k)
\]
lies in the $\Hmr\ull{\FF}_2$-Hurewicz image of $C(a^k)$ whenever $\lambda(k-i) \mid n$. In $\Hmr\ull{\FF}_2$-homology, \Cref{lem:iota_p_homology}(1) shows that $a^i u^n$ is supported on the bottom cell of $C(a^k)$, arising as the image of \[ a^i u^n \in \MM_{\star}\] under the bottom-cell inclusion $\upiota_{k\ast}$. However, by \Cref{main2}(1), its homological preimage $a^i u^n \in \MM_{\star}$ ($n > 0$) does \emph{not} lie in the $\Hmr\ull{\FF}_2$-Hurewicz image of $\SS_{\Cmr_2}$. This discrepancy naturally raises the question: how are the Hurewicz lifts of $a^i u^n$ in $\pi_{\star}^{\Cmr_2}(C(a^k))$ accounted for by the homotopy of $\SS_{\Cmr_2}$?

To answer this question, we examine the commutative diagram of exact sequences induced by the defining cofiber sequence \eqref{eq:full_cofiber_k} and the Hurewicz map $\Hfrak$:
\begin{equation} \label{eq:hurewicz_ladder}
\begin{tikzcd}[column sep=small]
 & & \textcolor{DarkBlue}{\widetilde{\upzeta}_{i,n}^{(k)}} \rar[color=DarkBlue, mapsto] \ar[color=DarkBlue, mapsto, ddd, bend left=11em] & \textcolor{DarkBlue}{{\sf p}_{k\ast}(\widetilde{\upzeta}_{i,n}^{(k)})} \ar[color=DarkBlue, mapsto, ddd, bend left=11em]\\
\dots \rar & \pi_{n-(n+i)\sigma}^{\Cmr_2}(\SS_{\Cmr_2}) \rar["\upiota_{k\ast}"] \dar["\Hfrak"'] & \pi_{n-(n+i)\sigma}^{\Cmr_2}(C(a^k)) \rar["{\sf p}_{k\ast}"] \dar["\Hfrak"] & \pi_{(n-1)-(n+i-k)\sigma}^{\Cmr_2}(\SS_{\Cmr_2}) \rar \dar["\Hfrak"] & \dots \\
\dots \rar & \MM_{n-(n+i)\sigma} \rar["\upiota_{k\ast}"'] & \Hmr_{n-(n+i)\sigma}^{\Cmr_2}(C(a^k); \ull{\FF}_2) \rar["{\sf p}_{k\ast}"'] & \MM_{(n-1)-(n+i-k)\sigma} \rar & \dots \\
 & \textcolor{DarkBlue}{a^i u^n} \rar[color=DarkBlue, mapsto] & \textcolor{DarkBlue}{a^i u^n} \rar[color=DarkBlue, mapsto] & \textcolor{DarkBlue}{0}
\end{tikzcd}
\end{equation}

\begin{notn}
Suppose $0 \le i < k$ and $n > 0$ satisfy $\lambda(k-i) \mid n$, so that $ a^i u^n \in \Hmr_{n-(n+i)\sigma}^{\Cmr_2}(C(a^k); \ull{\FF}_2)$ admits a Hurewicz lift \[ \widetilde{\upzeta}_{i,n}^{(k)} \in \pi_{n-(n+i)\sigma}^{\Cmr_2}(C(a^k))\] by \Cref{main1}. 
\end{notn}
\begin{prop} ${\sf p}_{k\ast}(\widetilde{\upzeta}_{i,n}^{(k)}) \in \pi_{(n-1)-(n+i-k)\sigma}^{\Cmr_2}(\SS_{\Cmr_2})$ is a nonzero class with trivial $\Hmr\ull{\FF}_2$-Hurewicz image. 
\end{prop}
\begin{proof} Since the top row in \eqref{eq:hurewicz_ladder} is exact, it is enough to show that $\widetilde{\upzeta}_{i,n}^{(k)}$ cannot lift along $\upiota_{k\ast}$ to $\pi_{n-(n+i)\sigma}^{\Cmr_2}(\SS_{\Cmr_2})$. Indeed, by \Cref{lem:iota_p_homology}(1), the homological map \[ \upiota_{k\ast} \colon \MM_{n-(n+i)\sigma} \to \Hmr_{n-(n+i)\sigma}^{\Cmr_2}(C(a^k); \ull{\FF}_2)\] in \eqref{eq:hurewicz_ladder} is an isomorphism in this degree. If $\widetilde{\upzeta}_{i,n}^{(k)} = \upiota_{k\ast}(x)$ for some $x \in \pi_{n-(n+i)\sigma}^{\Cmr_2}(\SS_{\Cmr_2})$, then commutativity of the left-hand square in \eqref{eq:hurewicz_ladder} would force $\Hfrak(x) = a^i u^n \in \MM_{n-(n+i)\sigma}$, contradicting \Cref{main2}(1). Thus
\[ 
{\sf p}_{k\ast}(\widetilde{\upzeta}_{i,n}^{(k)}) \neq 0 \in \pi_{(n-1)-(n+i-k)\sigma}^{\Cmr_2}(\SS_{\Cmr_2})
\]
is a non-trivial homotopy class in $\pi_{(n-1)-(n+i-k)\sigma}^{\Cmr_2}(\SS_{\Cmr_2})$. Moreover, the exactness of the bottom row of \eqref{eq:hurewicz_ladder} (also see \Cref{lem:iota_p_homology} (2)) implies  $\Hmr\ull{\FF}_2$-Hurewicz image of ${\sf p}_{k\ast}(\widetilde{\upzeta}_{i,n}^{(k)})$  
\begin{eqnarray*}
\mathfrak{H}({\sf p}_{k\ast}(\widetilde{\upzeta}_{i,n}^{(k)})) &=& {\sf p}_{k\ast}(\mathfrak{H}(\widetilde{\upzeta}_{i,n}^{(k)})) \\
&=& {\sf p}_{k\ast}(a^i u^n) \\ 
&=& 0
\end{eqnarray*}
is trivial. 
\end{proof}
\subsection{Relationships among the Hurewicz lifts $\widetilde{\upzeta}_{i,n}^{(k)}$} \ 

As $i$ and $k$ vary, the Hurewicz lifts $\widetilde{\upzeta}_{i,n}^{(k)}$ and their pinched images ${\sf p}_{k\ast}(\widetilde{\upzeta}_{i,n}^{(k)})$ are linked by the maps $\alpha^{(k)}$ of \eqref{eq:alpha} and $\beta^{(k)}$ of \eqref{eq:beta}. Extending the defining diagrams \eqref{eq:alpha} and \eqref{eq:beta} one step to the right yields the commutative diagrams of horizontal cofiber sequences
\begin{equation} \label{eq:alpha_pinch}
\begin{tikzcd}
\Sigma^{-(k+1)\sigma}\SS_{\Cmr_2} \rar["a^k"] \dar[equal] & \Sigma^{-\sigma}\SS_{\Cmr_2} \rar["\Sigma^{-\sigma}\upiota_k"] \dar["a"] & \Sigma^{-\sigma}C(a^k) \rar["\Sigma^{-\sigma}{\sf p}_k"] \dar["\alpha^{(k)}"] & \Sigma^{1-(k+1)\sigma}\SS_{\Cmr_2} \dar[equal] \\
\Sigma^{-(k+1)\sigma}\SS_{\Cmr_2} \rar["a^{k+1}"'] & \SS_{\Cmr_2} \rar["\upiota_{k+1}"'] & C(a^{k+1}) \rar["{\sf p}_{k+1}"'] & \Sigma^{1-(k+1)\sigma}\SS_{\Cmr_2}
\end{tikzcd}
\end{equation}
and
\begin{equation} \label{eq:beta_pinch}
\begin{tikzcd}
\Sigma^{-(k+1)\sigma}\SS_{\Cmr_2} \rar["a^{k+1}"] \dar["a"'] & \SS_{\Cmr_2} \rar["\upiota_{k+1}"] \dar[equal] & C(a^{k+1}) \rar["{\sf p}_{k+1}"] \dar["\beta^{(k)}"] & \Sigma^{1-(k+1)\sigma}\SS_{\Cmr_2} \dar["\Sigma^{1-k\sigma}a"] \\
\Sigma^{-k\sigma}\SS_{\Cmr_2} \rar["a^k"'] & \SS_{\Cmr_2} \rar["\upiota_k"'] & C(a^k) \rar["{\sf p}_k"'] & \Sigma^{1-k\sigma}\SS_{\Cmr_2}.
\end{tikzcd}
\end{equation}
In particular, the right-hand squares of \eqref{eq:alpha_pinch} and \eqref{eq:beta_pinch} give the relations
\begin{equation} \label{eq:p_alpha_beta_relations}
{\sf p}_{k+1} \circ \alpha^{(k)} \simeq \Sigma^{-\sigma}{\sf p}_k
\qquad \text{and} \qquad
{\sf p}_k \circ \beta^{(k)} \simeq a \circ {\sf p}_{k+1},
\end{equation}
while the middle squares imply that both composites $\alpha^{(k)} \circ \Sigma^{-\sigma}\beta^{(k)}$ and $\beta^{(k)} \circ \alpha^{(k)}$ are homotopic to multiplication by $a$.

A Hurewicz lift $\widetilde{\upzeta}_{i,n}^{(k)}$ of $a^i u^n \in \Hmr_{\star}^{\Cmr_2}(C(a^k); \ull{\FF}_2)$ is only well-defined modulo $\ker(\Hfrak)$. Therefore, we must take care when comparing lifts across different values of $i$ and $k$. The next two lemmas record how $\alpha^{(k)}$ and $\beta^{(k)}$ transport choices of Hurewicz lifts.

\begin{notn} \label{notn:zeta_standing}
Throughout this subsection, we fix $n, k > 0$ and $0 \le i < k$. We assume $\lambda(k-i) \mid n$ whenever considering $\widetilde{\upzeta}_{i,n}^{(k)} \in \pi_{n-(n+i)\sigma}^{\Cmr_2}(C(a^k))$, and $\lambda(k+1-i) \mid n$ whenever considering $\widetilde{\upzeta}_{i,n}^{(k+1)} \in \pi_{n-(n+i)\sigma}^{\Cmr_2}(C(a^{k+1}))$, so that the respective Hurewicz lifts of $a^i u^n$ exist by \Cref{main1}.
\end{notn}

\begin{lem}[{Invariance under $\alpha^{(k)}$}] \label{lem:zeta_alpha}
For any Hurewicz lift $\widetilde{\upzeta}_{i,n}^{(k)}$ of $a^i u^n$ in $C(a^k)$:
\begin{enumerate}
\item $\alpha^{(k)}_{\ast}(\widetilde{\upzeta}_{i,n}^{(k)}) \in \pi_{n-(n+i+1)\sigma}^{\Cmr_2}(C(a^{k+1}))$ is a Hurewicz lift of $a^{i+1} u^n$ in $C(a^{k+1})$.
\item ${\sf p}_{k+1\ast}\!\left(\alpha^{(k)}_{\ast}(\widetilde{\upzeta}_{i,n}^{(k)})\right) = {\sf p}_{k\ast}(\widetilde{\upzeta}_{i,n}^{(k)})$ in $\pi_{(n-1)-(n+i-k)\sigma}^{\Cmr_2}(\SS_{\Cmr_2})$.
\item Every Hurewicz lift $\widetilde{\upzeta}_{i+1,n}^{(k+1)}$ of $a^{i+1} u^n$ in $C(a^{k+1})$ is of the form $\alpha^{(k)}_{\ast}(\widetilde{\upzeta}_{i,n}^{(k)})$ for some lift $\widetilde{\upzeta}_{i,n}^{(k)}$ of $a^i u^n$ in $C(a^k)$.
\end{enumerate}
\end{lem}

\begin{proof}
Set $V = n - (n+i+1)\sigma$. By \eqref{eq:alpha_pinch} and naturality of $\Hfrak$, we have a commutative diagram
\[
\begin{tikzcd}[column sep=large]
\Hmr_{V+\sigma}^{\Cmr_2}(C(a^k); \ull{\FF}_2) \dar["\alpha^{(k)}_{\ast}"', "\cong"] & \pi_{V+\sigma}^{\Cmr_2}(C(a^k)) \lar["\Hfrak"'] \rar["{\sf p}_{k\ast}"] \dar[->>, "\alpha^{(k)}_{\ast}"] & \pi_{V + 1 - (k+1)\sigma}^{\Cmr_2}(\SS_{\Cmr_2}) \dar[equal] \\
\Hmr_{V}^{\Cmr_2}(C(a^{k+1}); \ull{\FF}_2) & \pi_{V}^{\Cmr_2}(C(a^{k+1})) \lar["\Hfrak"] \rar["{\sf p}_{k+1\ast}"'] & \pi_{V + 1 - (k+1)\sigma}^{\Cmr_2}(\SS_{\Cmr_2}).
\end{tikzcd}
\]
By \Cref{lem:iota_p_homology}(4), the left vertical map sends $a^i u^n \mapsto a^{i+1} u^n$ isomorphically, so the left square gives~(1) and the right square gives~(2). Since $|V| = -(i+1) < 0$, the middle vertical map $\alpha^{(k)}_{\ast}$ is surjective as in the proof of \Cref{main1}, yielding~(3).
\end{proof}

\begin{lem}[{Compatibility under $\beta^{(k)}$}] \label{lem:zeta_beta}
For any Hurewicz lift $\widetilde{\upzeta}_{i,n}^{(k+1)}$ of $a^i u^n$ in $C(a^{k+1})$:
\begin{enumerate}
\item $\beta^{(k)}_{\ast}(\widetilde{\upzeta}_{i,n}^{(k+1)}) \in \pi_{n-(n+i)\sigma}^{\Cmr_2}(C(a^k))$ is a Hurewicz lift of $a^i u^n$ in $C(a^k)$.
\item ${\sf p}_{k\ast}\!\left(\beta^{(k)}_{\ast}(\widetilde{\upzeta}_{i,n}^{(k+1)})\right) = a \cdot {\sf p}_{k+1\ast}(\widetilde{\upzeta}_{i,n}^{(k+1)})$.
\item $\alpha^{(k)}_{\ast}\!\left(\beta^{(k)}_{\ast}(\widetilde{\upzeta}_{i,n}^{(k+1)})\right) = a \cdot \widetilde{\upzeta}_{i,n}^{(k+1)}$, which is a Hurewicz lift of $a^{i+1} u^n$ in $C(a^{k+1})$.
\end{enumerate}
\end{lem}

\begin{proof}
Set $W = n - (n+i)\sigma$. By \eqref{eq:beta_pinch} and naturality of $\Hfrak$, we have a commutative diagram
\[
\begin{tikzcd}[column sep=large]
\Hmr_{W}^{\Cmr_2}(C(a^{k+1}); \ull{\FF}_2) \dar["\beta^{(k)}_{\ast}"'] & \pi_{W}^{\Cmr_2}(C(a^{k+1})) \lar["\Hfrak"'] \rar["{\sf p}_{k+1\ast}"] \dar["\beta^{(k)}_{\ast}"] & \pi_{W - 1 + (k+1)\sigma}^{\Cmr_2}(\SS_{\Cmr_2}) \dar["\cdot a"] \\
\Hmr_{W}^{\Cmr_2}(C(a^k); \ull{\FF}_2) & \pi_{W}^{\Cmr_2}(C(a^k)) \lar["\Hfrak"] \rar["{\sf p}_{k\ast}"'] & \pi_{W - 1 + k\sigma}^{\Cmr_2}(\SS_{\Cmr_2}).
\end{tikzcd}
\]
Since $0 \le i < k$, \Cref{lem:iota_p_homology}(3) shows that the left vertical map sends $a^i u^n \mapsto a^i u^n$. The left and right squares then give~(1) and~(2), while~(3) follows from \Cref{lem:zeta_alpha}(1) and the homotopy $\alpha^{(k)} \circ \Sigma^{-\sigma}\beta^{(k)} \simeq a \cdot \mr{id}_{C(a^{k+1})}$.
\end{proof}

\begin{rmk}[Coherent choices via maximal $a$-divisibility] \label{rmk:coherent_lifts}
Unlike $\alpha^{(k)}_{\ast}$, the map $\beta^{(k)}_{\ast}$ in the proof of \Cref{lem:zeta_beta} is neither injective nor surjective in general, as its fiber $\Sigma^{-k\sigma}C(a)$ and cofiber $\Sigma^{1-k\sigma}C(a)$ have homotopy groups $\pi_{k-i}^{e}(\SS)$ and $\pi_{k-1-i}^{e}(\SS)$ in non-negative stems. Consequently, an arbitrary Hurewicz lift $\widetilde{\upzeta}_{i,n}^{(k)}$ in $C(a^k)$ need not lie in the image of $\beta^{(k)}_{\ast}$.

To make our choices of $\widetilde{\upzeta}_{i,n}^{(k)}$ compatible across all $i$ and $k$, we exploit the equivalence 
\[ 
\lambda(k-i) \mid n \iff k - i \le \ell(n-1) + 1 
\]
from \Cref{notn:ell_lambda}. For each fixed $n > 0$, we first set $i = 0$ and $k = \ell(n-1) + 1$ (achieving the maximal gap $k - i = \ell(n-1) + 1$ allowed by \Cref{main1}), and fix a Hurewicz lift of $u^n \in \Hmr_{n(1-\sigma)}^{\Cmr_2}(C(a^{\ell(n-1)+1}); \ull{\FF}_2)$:
\[
\widetilde{\upzeta}_{0,n}^{(\ell(n-1)+1)} \in \pi_{n(1-\sigma)}^{\Cmr_2}\!\left(C(a^{\ell(n-1)+1})\right).
\]
For all $0 \le i < k$ satisfying $1 \le k - i \le \ell(n-1) + 1$, we then define $\widetilde{\upzeta}_{i,n}^{(k)}$ coherently by applying $\beta$ to decrease $k - i$ and $\alpha$ to shift $i$ and $k$:
\[
\widetilde{\upzeta}_{i,n}^{(k)} := \left(\alpha^{(k-1)}_{\ast} \circ \dots \circ \alpha^{(k-i)}_{\ast}\right) \circ \left(\beta^{(k-i)}_{\ast} \circ \dots \circ \beta^{(\ell(n-1))}_{\ast}\right)\!\left(\widetilde{\upzeta}_{0,n}^{(\ell(n-1)+1)}\right).
\]
By \Cref{lem:zeta_alpha} and \Cref{lem:zeta_beta}, these coherent lifts satisfy $\alpha^{(k)}_{\ast}(\widetilde{\upzeta}_{i,n}^{(k)}) = \widetilde{\upzeta}_{i+1,n}^{(k+1)}$ and $\beta^{(k)}_{\ast}(\widetilde{\upzeta}_{i,n}^{(k+1)}) = \widetilde{\upzeta}_{i,n}^{(k)}$ on the nose.
\end{rmk}

\subsection{Proof of \Cref{main3}} \ 

We now use the coherent Hurewicz lifts of \Cref{rmk:coherent_lifts} to construct the infinite family of elements $\zeta_{n-1}$ and establish our final main theorem:

\begin{thm*}[{\Cref{main3}}]
For every $n > 0$, there exists a non-zero element $\zeta_{n-1} \in \pi_{(n-1)(1 - \sigma)}^{\Cmr_2}(\SS_{\Cmr_2})$ such that:
\begin{enumerate}
\item $\zeta_{n-1}$ is $a$-torsion, i.e., $a \cdot \zeta_{n-1} = 0$ for all $n$.
\item $\zeta_{n-1}$ is divisible by $a^k$ if and only if $S^{n-1}$ admits $k$ linearly independent vector fields (i.e., $0 \le k \le \ell(n-1)$).
\item For any $k \ge 0$ such that $a^k$ divides $\zeta_{n-1}$, any choice of quotient $\frac{\zeta_{n-1}}{a^k}$ has a trivial $\Hmr\ull{\FF}_2$-Hurewicz image.
\end{enumerate}
\end{thm*}

\begin{proof}
Fix $n > 0$, and let $\ell := \ell(n-1)$ denote the maximal number of linearly independent vector fields on $S^{n-1}$. Let $\widetilde{\upzeta}_{0,n}^{(\ell+1)} \in \pi_{n(1-\sigma)}^{\Cmr_2}(C(a^{\ell+1}))$ be the Hurewicz lift of $u^n$ chosen in \Cref{rmk:coherent_lifts}, and define its pinched image
\[
\overline{\zeta}_{n-1} := {\sf p}_{\ell+1\ast}\!\left(\widetilde{\upzeta}_{0,n}^{(\ell+1)}\right) \in \pi_{(n-1)(1-\sigma) + \ell\sigma}^{\Cmr_2}(\SS_{\Cmr_2}).
\]
For each $0 \le k \le \ell$, the coherent lift $\widetilde{\upzeta}_{0,n}^{(k+1)} = (\beta^{(k+1)}_{\ast} \circ \dots \circ \beta^{(\ell)}_{\ast})(\widetilde{\upzeta}_{0,n}^{(\ell+1)}) \in \pi_{n(1-\sigma)}^{\Cmr_2}(C(a^{k+1}))$ of $u^n$ satisfies
\begin{equation} \label{eq:coherent_pinch}
{\sf p}_{k+1\ast}\!\left(\widetilde{\upzeta}_{0,n}^{(k+1)}\right) = a^{\ell - k} \cdot \overline{\zeta}_{n-1} \in \pi_{(n-1)(1-\sigma) + k\sigma}^{\Cmr_2}(\SS_{\Cmr_2})
\end{equation}
by \Cref{lem:zeta_beta}(2). Specializing \eqref{eq:coherent_pinch} to $k = 0$ (or equivalently, to any $\widetilde{\upzeta}_{i,n}^{(i+1)}$ via \Cref{lem:zeta_alpha}(2)), we define $\zeta_{n-1}$ to be the maximal non-zero $a$-multiple of $\overline{\zeta}_{n-1}$:
\[
\zeta_{n-1} := a^{\ell} \cdot \overline{\zeta}_{n-1} = {\sf p}_{1\ast}\!\left(\widetilde{\upzeta}_{0,n}^{(1)}\right) \in \pi_{(n-1)(1-\sigma)}^{\Cmr_2}(\SS_{\Cmr_2}).
\]

\medskip
\noindent
{\bf Proof of (1):} Since $\zeta_{n-1} = {\sf p}_{1\ast}(\widetilde{\upzeta}_{0,n}^{(1)})$, Proposition~4.2 implies that $\zeta_{n-1} \neq 0$. Moreover, by exactness of the defining cofiber sequence \eqref{eq:full_cofiber_k} for $C(a)$, the composite $a \circ {\sf p}_1$ is nullhomotopic, so $a \cdot \zeta_{n-1} = 0$.

\medskip
\noindent
{\bf Proof of (2):} If $S^{n-1}$ admits $k$ linearly independent vector fields ($0 \le k \le \ell$), then by \eqref{eq:coherent_pinch} the element $y_k := {\sf p}_{k+1\ast}(\widetilde{\upzeta}_{0,n}^{(k+1)}) = a^{\ell - k} \cdot \overline{\zeta}_{n-1}$ satisfies
\[
a^k \cdot y_k = a^{\ell} \cdot \overline{\zeta}_{n-1} = \zeta_{n-1},
\]
so $\zeta_{n-1}$ is divisible by $a^k$.

Conversely, suppose $\zeta_{n-1}$ is divisible by $a^k$ for some $k \ge 0$, so that $\zeta_{n-1} = a^k \cdot y$ for some $y \in \pi_{(n-1)(1-\sigma) + k\sigma}^{\Cmr_2}(\SS_{\Cmr_2})$. By~(1), we have
\[
a^{k+1} \cdot y = a \cdot (a^k \cdot y) = a \cdot \zeta_{n-1} = 0.
\]
By exactness of the cofiber sequence \eqref{eq:full_cofiber_k} for $C(a^{k+1})$, the class $y$ lifts along ${\sf p}_{k+1}$ to an element $\tilde{y} \in \pi_{n(1-\sigma)}^{\Cmr_2}(C(a^{k+1}))$ satisfying ${\sf p}_{k+1\ast}(\tilde{y}) = y$. Let $\beta^{(1,k)} := \beta^{(1)} \circ \dots \circ \beta^{(k)} \colon C(a^{k+1}) \to C(a)$ denote the iterated reduction map. Since ${\sf p}_1 \circ \beta^{(1,k)} \simeq a^k \circ {\sf p}_{k+1}$, we have
\[
{\sf p}_{1\ast}\!\left(\beta^{(1,k)}_{\ast}(\tilde{y})\right) = a^k \cdot {\sf p}_{k+1\ast}(\tilde{y}) = a^k \cdot y = \zeta_{n-1} = {\sf p}_{1\ast}\!\left(\widetilde{\upzeta}_{0,n}^{(1)}\right).
\]
By exactness of \eqref{eq:hurewicz_ladder} for $C(a)$, the difference $\beta^{(1,k)}_{\ast}(\tilde{y}) - \widetilde{\upzeta}_{0,n}^{(1)}$ lies in $\ker({\sf p}_{1\ast}) = \mr{Im}(\upiota_{1\ast})$, so $\beta^{(1,k)}_{\ast}(\tilde{y}) = \widetilde{\upzeta}_{0,n}^{(1)} + \upiota_{1\ast}(x)$ for some $x \in \pi_{n(1-\sigma)}^{\Cmr_2}(\SS_{\Cmr_2})$. Because $n > 0$, \Cref{main2}(1) implies that $\Hfrak(x) = 0 \in \MM_{n(1-\sigma)}$. Applying $\Hfrak$ therefore gives
\[
\beta^{(1,k)}_{\ast}(\Hfrak(\tilde{y})) = \Hfrak\!\left(\beta^{(1,k)}_{\ast}(\tilde{y})\right) = \Hfrak\!\left(\widetilde{\upzeta}_{0,n}^{(1)}\right) + \upiota_{1\ast}(\Hfrak(x)) = u^n \in \Hmr_{n(1-\sigma)}^{\Cmr_2}(C(a); \ull{\FF}_2).
\]
Since $\beta^{(1,k)}_{\ast} \colon \Hmr_{n(1-\sigma)}^{\Cmr_2}(C(a^{k+1}); \ull{\FF}_2) \to \Hmr_{n(1-\sigma)}^{\Cmr_2}(C(a); \ull{\FF}_2)$ is an isomorphism (both groups being $\FF_2\{u^n\}$), it follows that $\Hfrak(\tilde{y}) = u^n$. Thus $u^n$ lies in the $\Hmr\ull{\FF}_2$-Hurewicz image of $C(a^{k+1})$, which by \Cref{main1} holds if and only if $S^{n-1}$ admits $k$ linearly independent vector fields ($0 \le k \le \ell $).

\medskip
\noindent
{\bf Proof of (3):} Let $y \in \pi_{(n-1)(1-\sigma) + k\sigma}^{\Cmr_2}(\SS_{\Cmr_2})$ be any choice of quotient $\frac{\zeta_{n-1}}{a^k}$, so that $a^k \cdot y = \zeta_{n-1}$. As shown in the proof of~(2), $y = {\sf p}_{k+1\ast}(\tilde{y})$ for a Hurewicz lift $\tilde{y} \in \pi_{n(1-\sigma)}^{\Cmr_2}(C(a^{k+1}))$ of $u^n$. By Proposition~4.2 (or \Cref{lem:iota_p_homology}(2)), we conclude that
\[
\Hfrak(y) = \Hfrak({\sf p}_{k+1\ast}(\tilde{y})) = {\sf p}_{k+1\ast}(\Hfrak(\tilde{y})) = {\sf p}_{k+1\ast}(u^n) = 0. \qedhere
\]
\end{proof}

\bibliographystyle{amsalpha}
\bibliography{C2-Hurewicz}

@unpublished{GI,
  title={C2-equivariant stable stems},
  author={Guillou, Bertrand J. and  Isaksen, Daniel C.},
  journal={},
  note = {Available at \url{https://arxiv.org/abs/2404.14627} },
}

@article {BBQ,
    AUTHOR = {Bhattacharya, Prasit and Bobkova, Irina and Quigley, J. D.},
     TITLE = {New infinite families in the stable homotopy groups of
              spheres},
   JOURNAL = {Geom. Topol.},
  FJOURNAL = {Geometry \& Topology},
    VOLUME = {30},
      YEAR = {2026},
    NUMBER = {6},
     PAGES = {2367--2393},
      ISSN = {1465-3060,1364-0380},
   MRCLASS = {55Q10 (55Q45 55Q51)},
  MRNUMBER = {5111600},
       DOI = {10.2140/gt.2026.30.2367},
       URL = {https://doi-org.nmsu.idm.oclc.org/10.2140/gt.2026.30.2367},
}

@article {DuggerKR,
    AUTHOR = {Dugger, Daniel},
     TITLE = {An {A}tiyah-{H}irzebruch spectral sequence for {$KR$}-theory},
   JOURNAL = {$K$-Theory},
  FJOURNAL = {$K$-Theory. An Interdisciplinary Journal for the Development,
              Application, and Influence of $K$-Theory in the Mathematical
              Sciences},
    VOLUME = {35},
      YEAR = {2005},
    NUMBER = {3-4},
     PAGES = {213--256},
      ISSN = {0920-3036,1573-0514},
   MRCLASS = {19L64 (55S45 55T25)},
  MRNUMBER = {2240234},
MRREVIEWER = {G\'erald\ Gaudens},
       DOI = {10.1007/s10977-005-1552-9},
       URL = {https://doi-org.nmsu.idm.oclc.org/10.1007/s10977-005-1552-9},
}

@incollection {LewisCP,
    AUTHOR = {Lewis, Jr., L. Gaunce},
     TITLE = {The {$R{\rm O}(G)$}-graded equivariant ordinary cohomology of
              complex projective spaces with linear {${\bf Z}/p$} actions},
 BOOKTITLE = {Algebraic topology and transformation groups ({G}\"ottingen,
              1987)},
    SERIES = {Lecture Notes in Math.},
    VOLUME = {1361},
     PAGES = {53--122},
 PUBLISHER = {Springer, Berlin},
      YEAR = {1988},
      ISBN = {3-540-50528-8},
   MRCLASS = {55N91 (55N22)},
  MRNUMBER = {979507},
       DOI = {10.1007/BFb0083034},
       URL = {https://doi-org.nmsu.idm.oclc.org/10.1007/BFb0083034},
}

@article {HuKriz,
    AUTHOR = {Hu, Po and Kriz, Igor},
     TITLE = {Real-oriented homotopy theory and an analogue of the
              {A}dams-{N}ovikov spectral sequence},
   JOURNAL = {Topology},
  FJOURNAL = {Topology. An International Journal of Mathematics},
    VOLUME = {40},
      YEAR = {2001},
    NUMBER = {2},
     PAGES = {317--399},
      ISSN = {0040-9383},
   MRCLASS = {55T15 (19L47 55N22 55N91 55P42 55P91)},
  MRNUMBER = {1808224},
MRREVIEWER = {J.\ P. C. Greenlees},
       DOI = {10.1016/S0040-9383(99)00065-8},
       URL = {https://doi-org.nmsu.idm.oclc.org/10.1016/S0040-9383(99)00065-8},
}

@article {Clover,
    AUTHOR = {May, Clover},
     TITLE = {A structure theorem for {$RO(C_2)$}-graded {B}redon
              cohomology},
   JOURNAL = {Algebr. Geom. Topol.},
  FJOURNAL = {Algebraic \& Geometric Topology},
    VOLUME = {20},
      YEAR = {2020},
    NUMBER = {4},
     PAGES = {1691--1728},
      ISSN = {1472-2747,1472-2739},
   MRCLASS = {55N91},
  MRNUMBER = {4127082},
MRREVIEWER = {Nansen\ Petrosyan},
       DOI = {10.2140/agt.2020.20.1691},
       URL = {https://doi-org.nmsu.idm.oclc.org/10.2140/agt.2020.20.1691},
}

@article {BEM,
    AUTHOR = {Bhattacharya, Prasit and Egger, Philip and Mahowald, Mark},
     TITLE = {On the periodic {$v_2$}-self-map of {$A_1$}},
   JOURNAL = {Algebr. Geom. Topol.},
  FJOURNAL = {Algebraic \& Geometric Topology},
    VOLUME = {17},
      YEAR = {2017},
    NUMBER = {2},
     PAGES = {657--692},
      ISSN = {1472-2747,1472-2739},
   MRCLASS = {55Q51},
  MRNUMBER = {3623667},
MRREVIEWER = {Gijs\ Heuts},
       DOI = {10.2140/agt.2017.17.657},
       URL = {https://doi-org.nmsu.idm.oclc.org/10.2140/agt.2017.17.657},
}

@article {AdamsVF,
    AUTHOR = {Adams, J. F.},
     TITLE = {Vector fields on spheres},
   JOURNAL = {Ann. of Math. (2)},
  FJOURNAL = {Annals of Mathematics. Second Series},
    VOLUME = {75},
      YEAR = {1962},
     PAGES = {603--632},
      ISSN = {0003-486X},
   MRCLASS = {57.30},
  MRNUMBER = {139178},
MRREVIEWER = {M.\ F.\ Atiyah},
       DOI = {10.2307/1970213},
       URL = {https://doi-org.nmsu.idm.oclc.org/10.2307/1970213},
}

@unpublished{HurHF,
  title={A Hurewicz theorem for $RO(C_2)$-graded equivariant homology governed by vector fields of spheres},
  author={ Guo, Manyi and Li, Guchuan and Lu, Yunze and Ma, Sihao and Wu, Yuchen and Xu, Zhouli and Yang, Albert Jinghui and Zhang, Shangjie},
  journal={},
  note = {Available at \url{https://arxiv.org/abs/2605.26334} },
}

@unpublished{HurKR,
  title={The $RO(C_2)$-graded Hurewicz image of Atiyah's Real K-theory},
  author={  Ma, Sihao and Shi, Danny Xiaolin  and Xu, Zhouli and Yang, Albert Jinghui and Zhang, Shangjie},
  journal={},
  note = {Available at \url{https://arxiv.org/abs/2610.19606} },
}

@article {James,
    AUTHOR = {James, I. M.},
     TITLE = {Spaces associated with {S}tiefel manifolds},
   JOURNAL = {Proc. London Math. Soc. (3)},
  FJOURNAL = {Proceedings of the London Mathematical Society. Third Series},
    VOLUME = {9},
      YEAR = {1959},
     PAGES = {115--140},
      ISSN = {0024-6115,1460-244X},
   MRCLASS = {55.00},
  MRNUMBER = {102810},
MRREVIEWER = {W.\ S.\ Massey},
       DOI = {10.1112/plms/s3-9.1.115},
       URL = {https://doi-org.nmsu.idm.oclc.org/10.1112/plms/s3-9.1.115},
}

@unpublished{LWX,
  title={The $RO(C_2)$-graded Hurewicz image of Atiyah's Real K-theory},
  author={    Lin, Weinan and Wang, Guozhen and  Xu, Zhouli },
  journal={},
  note = {Available at \url{https://arxiv.org/abs/2412.10879} },
}

@article {Pham,
    AUTHOR = {Pham, Viet-Cuong},
     TITLE = {On the surjectivity of the tmf-{H}urewicz image of {$A_1$}},
   JOURNAL = {Algebr. Geom. Topol.},
  FJOURNAL = {Algebraic \& Geometric Topology},
    VOLUME = {23},
      YEAR = {2023},
    NUMBER = {1},
     PAGES = {217--241},
      ISSN = {1472-2747,1472-2739},
   MRCLASS = {55Q10},
  MRNUMBER = {4568003},
MRREVIEWER = {Prasit\ Bhattacharya},
       DOI = {10.2140/agt.2023.23.217},
       URL = {https://doi-org.nmsu.idm.oclc.org/10.2140/agt.2023.23.217},
}

@article {BehShah,
    AUTHOR = {Behrens, Mark and Shah, Jay},
     TITLE = {{$C_2$}-equivariant stable homotopy from real motivic stable
              homotopy},
   JOURNAL = {Ann. K-Theory},
  FJOURNAL = {Annals of K-Theory},
    VOLUME = {5},
      YEAR = {2020},
    NUMBER = {3},
     PAGES = {411--464},
      ISSN = {2379-1683,2379-1691},
   MRCLASS = {14F42 (55N91 55P91 55Q91)},
  MRNUMBER = {4132743},
MRREVIEWER = {Stephen\ McKean},
       DOI = {10.2140/akt.2020.5.411},
       URL = {https://doi-org.nmsu.idm.oclc.org/10.2140/akt.2020.5.411},
}

@unpublished{Ma,
  title={The Borel and genuine $C_2$-equivariant Adams spectral sequences},
  author={  Ma, Sihao},
  journal={},
  note = {Available at \url{https://arxiv.org/abs/2208.12883} },
}

@article {DI,
    AUTHOR = {Dugger, Daniel and Isaksen, Daniel C.},
     TITLE = {Low-dimensional {M}ilnor-{W}itt stems over {$\Bbb{R}$}},
   JOURNAL = {Ann. K-Theory},
  FJOURNAL = {Annals of K-Theory},
    VOLUME = {2},
      YEAR = {2017},
    NUMBER = {2},
     PAGES = {175--210},
      ISSN = {2379-1683,2379-1691},
   MRCLASS = {14F42 (55Q45 55S10 55T15)},
  MRNUMBER = {3590344},
MRREVIEWER = {Oliver\ R\"ondigs},
       DOI = {10.2140/akt.2017.2.175},
       URL = {https://doi-org.nmsu.idm.oclc.org/10.2140/akt.2017.2.175},
}

@article {AtiyahThom,
    AUTHOR = {Atiyah, M. F.},
     TITLE = {Thom complexes},
   JOURNAL = {Proc. London Math. Soc. (3)},
  FJOURNAL = {Proceedings of the London Mathematical Society. Third Series},
    VOLUME = {11},
      YEAR = {1961},
     PAGES = {291--310},
      ISSN = {0024-6115},
   MRCLASS = {57.30},
  MRNUMBER = {131880},
MRREVIEWER = {R. Bott},
       DOI = {10.1112/plms/s3-11.1.291},
       URL = {https://doi-org.proxy.library.nd.edu/10.1112/plms/s3-11.1.291},
}

\end{document}